\documentclass[10pt]{article}
\usepackage[latin1]{inputenc}
\usepackage{amscd}
\usepackage{amsfonts}
\usepackage{amsgen}
\usepackage{amsmath}
\usepackage{amssymb}
\usepackage{amstext}
\usepackage{amsthm}
\usepackage{graphicx}
\usepackage{indentfirst}
\usepackage{latexsym}
\usepackage{makeidx}
\usepackage{mathrsfs}
\usepackage{epsfig}
\usepackage[all]{xy}
\usepackage[usenames]{color}
\usepackage{enumerate}

\newlength{\dede}     

\newtheorem{Df}{Definition}[section]
\newtheorem{Teo}[Df]{Theorem}
\newtheorem{Prop}[Df]{Proposition}
\newtheorem{Lem}[Df]{Lemma}
\newtheorem{Ex}[Df]{Example}
\newtheorem{Obs}[Df]{Remark}

\newtheorem{Que}[Df]{Question}
\newtheorem{Cor}[Df]{Corollary}

\newcommand{\n}{\noindent}
\newcommand{\Dem}{\n{\bf Proof:\;\;}}

\newcommand{\bc}{\begin{center}}
\newcommand{\ec}{\end{center}}

\newcommand{\N}{\mathbb{N}}

\newcommand{\vtres}{\vspace*{0.3cm}}

\def \beq { \begin{equation} }
\def \eeq { \end{equation} }

\newcommand{\baf}{\begin{Afi}\nf{\sl . }}
\newcommand{\eaf}{\end{Afi}}

\newcommand{\bdf}{\begin{Df}\nf{\bf .}}
\newcommand{\edf}{\end{Df}}
\newcommand{\bte}{\begin{Th}\nf{\bf .}}
\newcommand{\ete}{\end{Th}}  
\newcommand{\bco}{\begin{Co}\nf{\bf .}}
\newcommand{\eco}{\end{Co}}
\newcommand{\ble}{\begin{Le}\nf{\bf .}}
\newcommand{\ele}{\end{Le}}
\newcommand{\bpr}{\begin{Pro}\nf{\bf .}}
\newcommand{\epr}{\end{Pro}}
\newcommand{\bex}{\begin{Exa}\nf{\bf .} \rm}
\newcommand{\eex}{\end{Exa}}
\newcommand{\bre}{\begin{Rem}\nf{\bf .} \rm}
\newcommand{\ere}{\end{Rem}}
\newcommand{\bfa}{\begin{Fa}\nf{\bf .} \sl}
\newcommand{\efa}{\end{Fa}}
\newcommand{\bqt}{\begin{Que}\nf{\bf .}}
\newcommand{\eqt}{\end{Que}}

\newcommand{\bct}{\begin{Ct}\nf{\bf .} \rm}
\newcommand{\ect}{\end{Ct}}
\newcommand{\bxa}{\begin{Exa}\nf{\bf .} \rm}
\newcommand{\exa}{\end{Exa}}

\newcommand{\Fm}[2]{Fm_{#1}(#2)}   

\begin{document}

\title{Congruence Filter Pairs, Adjoints and Leibniz Hierarchy}
\author{P. Arndt, 
H.L. Mariano, D.C. Pinto}

\maketitle




\begin{abstract}
We have introduced in \cite{AMP1} the notion of (finitary) filter pair as a tool for creating and analyzing logics. A  filter pair can be seen as a presentation of a logic, given by presenting its lattice of theories as the image of a lattice homomorphism, with certain properties ensuring that the resulting logic is finitary and substitution invariant.

Every finitary, substitution invariant logic arises from a filter pair. Particular classes of logics can be characterized as arising from special classes of filter pairs. We consider so-called congruence filter pairs, i.e. filter pairs for which the domain of the lattice homomorphism is a lattice of congruences for some quasivariety. We show that the class of logics admitting a presentation by such a filter pair is exactly the class of logics having an algebraic semantics.

We study the properties of a certain Galois connection coming with such filter pairs. 
We give criteria for a congruence filter pair to present a logic in some classes of the Leibniz hierarchy by means of this Galois connection, and its interplay with the Leibniz operator. 

As an application, we show a bridge theorem, stating that the amalgamation property implies the Craig interpolation property, for a certain class of logics including non-protoalgebraic logics.

\end{abstract}

\section*{Introduction}

The first main aim of this article is to introduce the notion of filter pair as a tool for creating and analyzing logics. A filter pair can be seen as a presentation of a logic, different in style from the usual syntactic presentations in terms of axioms and rules, or semantic presentations in terms of matrices. Roughly, a filter pair is a highly structured lattice homomorphism from some algebraic lattice into the power set of the formula algebra, whose image can be taken as the lattice of theories of a logic.
The notion of filter pair can be developed into several different directions, some of which are sketched in Section \ref{SectionVista}.

In this article we concentrate on one direction, which constitutes its second main aim: to advance the general study of logics having an algebraic semantics. Recall that an algebraic semantics for a logic $l$ is given by a class of algebras ${\bf K}$ and a set of equations $\delta(x) = \epsilon(x)$ such that $$\Gamma \vdash_l \varphi \ \ \ \ \Leftrightarrow \ \ \ \ \{\delta(\gamma) = \epsilon(\gamma) \mid \gamma \in \Gamma\} \vDash_{{\bf K}} \delta(\varphi) = \epsilon(\varphi).$$
The notion was considered by Blok-Pigozzi in \cite{BP1}, but there mainly used as a stepping stone towards the notion of algebraizable logic. 
As Font writes \cite[p.114]{Fon}: ``The bare concept of algebraic semantics has not been extensively studied. [...] At the time of writing, \cite{BlokRebagliato} and \cite{Raftery} are probably the only papers that contain significant results about it.''
It turns out that logics having an algebraic semantics are exactly those which can be presented by a so-called equational filter pair - in the above picture of a filter pair as a highly structured lattice homomorphism, these are those homomorphisms whose domain is the lattice of congruences relative to some quasivariety ${\bf K}$. We find that such a presentation, and the additional structure that comes with it, can be employed to analyze logics with an algebraic semantics. We illustrate this by considering some criteria for being algebraizable or truth-equational and for satisfying the Craig interpolation property.

\medskip

We now sketch the basic idea of the notion of filter pair:

Throughout the article the word \emph{logic} will mean a pair $(\Sigma, \vdash)$ where $\Sigma$ is a signature, i.e. a collection of connectives with finite arities, and $\vdash$ is a Tarskian consequence relation, i.e. an idempotent, increasing, monotonic, finitary and structural relation between subsets and elements of the set of formulas $\Fm{\Sigma}{X}$ built from $\Sigma$ and a set $X$ of variables.

It is well-known that every Tarskian logic gives rise to an algebraic lattice contained in the powerset $\wp(\Fm{\Sigma}{X})$, namely the lattice of theories. This lattice is closed under arbitrary intersections and directed unions. 

Conversely an algebraic lattice $L \subseteq \wp(\Fm{\Sigma}{X})$ that is closed under arbitrary intersections and unions of increasing chains gives rise to a finitary closure operator (assigning to a subset $A \subseteq \Fm{\Sigma}{X}$ the intersection of all members of $L$ containing $A$). This closure operator need not be structural --- this is an extra requirement.

We observe that the structurality of the logic just defined is equivalent to the \emph{naturality} (in the sense of category theory) of the inclusion of the algebraic lattice into the power set of formulas with respect to endomorphisms of the formula algebra: Structurality means that the preimage under a substitution of a theory is a theory again or, equivalently, that the following diagram commutes:

$$\xymatrix{
\Fm{\Sigma}{X} \ar[d]_\sigma & L \ar@{^(->}[r]^-i & \wp(\Fm{\Sigma}{X}) \\
\Fm{\Sigma}{X} & L \ar[u]^{\sigma^{-1}|_{_L}} \ar@{^(->}[r]^-i & \wp(\Fm{\Sigma}{X}) \ar[u]_{\sigma^{-1}}
}$$

Further, it is equivalent to demand this naturality for all $\Sigma$-algebras and homomorphisms instead of just the formula algebra.

We thus arrive at the definition of \emph{filter pair}, Def. \ref{DefinitionFilterPair}: A filter pair for the signature $\Sigma$ is a contravariant functor $G$ from $\Sigma$-algebras to algebraic lattices together with a natural transformation $i \colon G \to \wp(-)$ from $G$ to the functor taking an algebra to the power set of its underlying set, which preserves arbitrary infima and directed suprema.

The logic associated to a filter pair $(G, i)$ is simply the logic associated (in the above fashion) to the algebraic lattice given by the image $i(G(\Fm{\Sigma}{X})) \subseteq \wp(\Fm{\Sigma}{X})$.

In particular, it is clear that different filter pairs can give rise to the same logic, indeed this will happen precisely if the images of $i$ for the formula algebra are the same. A filter pair can thus be seen as a \emph{presentation} of a logic, and there can of course be different presentations of the same logic.
We could have removed a bit of this ambiguity by requiring that $i$ be an inclusion, but it is one of the insights of this article that it is beneficial not to do this. Indeed this will give us greater flexibility for the choice of the functor $G$, and injectivity of $i$ can become a meaningful extra feature. Thus, for example, if $G$ is the functor associating to a $\Sigma$-structure the lattice of congruences relative to some quasivariety ${\bf K}$, i.e. when we have a so-called \emph{congruence filter pair}, then by Prop. \ref{2.3} the injectivity of $i$ means that the associated logic is algebraizable.

In this article we show how to recognize classes of logics, like algebraizable or truth-equational logics, through their presentations by filter pairs and how these presentations may permit to use algebraic methods even outside the realm of protoalgebraic logics, see the following overview for a sample of results.

Another theme that we take up is the theme of adjunctions arising in the context of filter pairs. It follows from the definition that the maps $i$ in a filter pair each have a left adjoint $\Xi$ (i.e. are part of a Galois connection). The closure operator obtained as the composition $i \circ \Xi$ is exactly the closure operator generating the theories of the logic. For congruence filter pairs we also have another operator back from theories to congruences, namely the Leibniz operator $\Omega$. The interplay of the three operators $i, \Xi, \Omega$ gives an interesting picture that one can observe for any logic with an algebraic semantics. As an illustration of the concrete relevance of these operators, we show in Section \ref{SectionCraigInterpolation} that whenever the maps $\Xi$ form a natural transformation (actually a weaker condition suffices), then the amalgamation property of a quasivariety implies the Craig interpolation property of a logic with an algebraic semantics in this quasivariety. While we have not yet been able to make this result bear fruits for the most common types of logics, we can already use it to cover a certain class of examples which includes non-protoalgebraic and non-truth-equational logics, showing that the results on amalgamation and interpolation can be extended beyond the known cases from the literature -- see Def. \ref{DefregularEquationsAndVarieties} and Cor. \ref{CorollaryAmalgamationInregularVarietiesImpliesInterpolation}. 

\medskip

{\bf Overview of the article:} 

In section \ref{SectionPreliminaries} we fix terminology and gather some results characterizing different classes of logics, like algebraizable, truth-equational or protoalgebraic logics. No original results are contained in this section.

In section \ref{SectionFilterFunctors} we introduce the notion of filter pair (Def. \ref{DefinitionFilterPair}) and show how to get a logic from a filter pair (Prop. \ref{LogicsFromFilterPairs}). We show some additional structures that one obtains from a filter pair, for example a left adjoint to each map $i_A \ (A \in \Sigma\text{-Str})$ (Def. \ref{DefLeftAdjointOfAFilterPair}) whose composition with $i$ yields the closure operator of the logic (Prop. \ref{PropDescriptionOfConsequenceRelationByClosureOperator}) and a consequence relation on each $A \in \Sigma\text{-Str}$ (Prop. \ref{LogicsFromFilterPairs}). We establish the basic relations between these structures, e.g. that homomorphisms of $\Sigma$-structures induce translations between these consequence relations, which are conservative in the case of free algebras (Prop. \ref{PropHomomorphismsInduceTranslationsAndVariableInclusionsInduceConservativeTranslations}), that subsets in the image of $i$ are filters for the logic, and that in the case of free algebras all filters are in the image of $i$ (Prop. \ref{FreeAlgebraResult}).

In section \ref{SectionFilterPairsOverCo} we consider the special case of filter pairs where the functor $G$ is given by congruences relative to some class of algebras. In subsection 3.1 we consider the question what the properties of $G$ demanded in the definition of filter pair enforce on the class ${\bf K}$. We find that ${\bf K}$ must be closed under subalgebras (Prop. \ref{PropCoKIsAFunctor}) and products (Prop. \ref{PropKClosedUnderProducts}). This subsection may be skipped and ${\bf K}$ be assumed to be a quasivariety in the remainder of the article.
In subsection 3.2 we introduce congruence filter pairs, and give a general construction of them, yielding the so-called equational filter pairs (Thm. \ref{TheoremLogicsFromEquations}). We give criteria when the associated logic is algebraizable (Prop. \ref{CriterionAlgebraizable}) or truth-equational (Cor. \ref{CriterionTruthEquational}, Prop. \ref{CriterionTruthEquationalPointedQuasiVar}). 

In section \ref{filteradjunctionequationalfilterpairs} we show and analyze some of the extra structure coming with an equational filter pair, for example give more concrete descriptions of the left adjoint of the $i_A$ mentioned in the previous paragraph (Prop. \ref{PropFirstFormulaForLeftAdjointForEquationalFilterPairs}, Prop. \ref{PropOurLeftAdjointIsBlokPigozzisOmegaK}) and investigate its relation with the Leibniz operator (Prop. \ref{TheoremInclusionsXiIdLeibniz} and the following discussion). We conclude from these results that logics having an algebraic semantics are precisely those which admit a presentation by an equational filter pair (Thm. \ref{TeoEveryLogicWithAlgebraicSemanticsComesFromEqunlFilterPair}), and we get another criterion for the algebraizability of a logic (Thm. \ref{TheoremInclusionsXiIdLeibniz}).

In section \ref{SectionCraigInterpolation} we show how the language of filter pairs can be used to establish connections between the amalgamation property of a quasivariety and Craig interpolation of a logic with algebraic semantics in that quasivariety, Theorem \ref{TheoremMatrixAmalgamationImpliesCraigInterpolation}. We recover the fact that amalgamation implies interpolation for algebraizable logics (Cor. \ref{CorollaryForAlgebraizableLogicsAmalgamationImpliesCraigInterpolation}), but also cover some cases which include  non-protoalgebraic logics (Lem. \ref{LemmaConditionForNaturalityOfXi}, Cor. \ref{CorollaryAmalgamationInregularVarietiesImpliesInterpolation}).

Finally, in section \ref{SectionVista} we discuss some of the other directions into which one can develop filter pairs.





\section{Preliminaries}\label{SectionPreliminaries}

\subsection{Ordered sets}

\begin{Df}[Galois connections]
 Let $f \colon P \to Q$, $g \colon Q \to P$ be order preserving maps between partially ordered sets. Then $(f,g)$ is a Galois connection, where  $f$ is called left adjoint to $g$ and $g$ is called right adjoint to $f$, if the following equivalence is satisfied: 
\[f(p) \leq q \Leftrightarrow  p \leq g(q) \ \ \ \ \forall p \in P,\ q \in Q\]
\end{Df}

A  Galois connection (or adjunction) between posets, $(f,g)$, is uniquely determined by its left (respectively right) component,  and  induces an isomorphism between the subposets  $\hat{P} \subseteq P$ and $\hat{Q} \subseteq  Q$, where $\hat{P} := Fix(g \circ f)$ and $\hat{Q} := Fix(f \circ g)$.

\begin{Teo}\label{TheoremExistenceOfLeftAdjoint}
{\em \cite[Thm. 3.6.9]{TaylorPracticalFoundations}} Let $g \colon Q \to P$ be a function between complete (and hence cocomplete) partially ordered sets that preserves arbitrary infima (resp. suprema). Then $g$ has a left adjoint (resp. right adjoint) $f \colon P \to Q$, given by $f(p)=inf\{q \in Q \, \mid \, p \leq g(q) \}$ (resp. $f(p) = sup\{q \in Q \, \mid \, p \geq g(q) \}$).
\end{Teo}

We will be concerned with algebraic lattices.

\begin{Df}[algebraic lattices] \label{AL-df}
\begin{enumerate}
\item Let $L$ be a lattice. A element $a\in L$ is compact if for every (upward) directed subset $\{d_{i}\}$ of $L$ we have $a\leq \bigvee_{i}d_{i}\ \Leftrightarrow\ \exists i(a\leq d_{i})$. $L$ is said to be  algebraic if it is a complete lattice such that every element is the join of the compact elements below it. 
\item We denote by $AL$ the category of all algebraic lattices and morphisms that preserve arbitrary meets and directed joins. In particular, these morphisms are monotonic functions that preserve the top elements of the lattices.
\end{enumerate}
\end{Df}

\begin{Obs}[sublattices of algebraic lattices] \label{AL-obs}
Let $X$ be a set and denote $(P(X), \subseteq)$ the (algebraic) lattice of all subsets of $X$. Note that a subset ${\cal S} \subseteq P(X)$ that is closed under arbitrary intersections and directed unions determines an algebraic sublattice $({\cal S}, \subseteq) \hookrightarrow (P(X), \subseteq)$, where the compact elements coincide with the finitely generated elements of $({\cal S}, \subseteq)$.
\end{Obs}

\subsection{Structures, logics, filters and matrices}

\begin{Df}[signatures, structures, formulas]
A signature is a sequence of pairwise disjoint sets $\Sigma=(\Sigma_{n})_{n\in \N}$, where $\Sigma_{n}$ is thought of as a set of $n$-ary function symbols. A $\Sigma$-structure is a set with an interpretation of the function symbols as actual operations. A homomorphism of $\Sigma$-structures is a map respecting these operations. We denote by $\Sigma$-Str the category of $\Sigma$-structures and homomorphisms.

For a set $X$ denote by $\Fm{\Sigma}{X}$ (or $\Fm{}{X}$, if the signature is clear), the absolutely free $\Sigma$-structure, i.e. the set of $\Sigma$-formulas with variables from $X$. A homomorphism $\Fm{}{X} \to A$ with domain $\Fm{}{X}$ is also called a valuation. An endomorphism of $\Fm{\Sigma}{X}$ is also called a substitution.
\end{Df}

In this work there will be naturally occurring consequence relations on non-free $\Sigma$-structures, hence we phrase the definition of logic slightly more generally than customary.

\begin{Df}[consequence relation, abstract logic, logic]\label{DefAbstractLogic}
A Tarskian consequence relation on a $\Sigma$-structure $A$ is a relation $\vdash\subseteq\wp(A)\times A$, such that, for all sets $\Gamma,\Delta \subseteq A$ and all elements $\varphi,\psi \in A$ the following conditions are satisfied:
\begin{itemize}
\item[$\circ$]$\bf{Reflexivity}: $If $\varphi\in\Gamma,\ \Gamma\vdash\varphi$
\item[$\circ$]$\bf{Cut}: $If $\Gamma\vdash\varphi$ and for every $\psi\in\Gamma,\ \Delta\vdash\psi$, then $\Delta\vdash\varphi$
\item[$\circ$]$\bf{Monotonicity}: $If $\Gamma\subseteq\Delta$ and $\Gamma\vdash\varphi$, then $\Delta\vdash\varphi$
\item[$\circ$]$\bf{Finitarity}: $If $\Gamma\vdash\varphi$, then there is a finite subset $\Delta$ of $\Gamma$ such that $\Delta\vdash\varphi$.
\item[$\circ$]$\bf{Structurality}: $If $\Gamma\vdash\varphi$ and $\sigma$ is an endomorphism of $A$, then $\sigma[\Gamma]\vdash\sigma(\varphi)$
\end{itemize}

An abstract logic (in the sense of \cite{BloomBrownSuszko}) is a pair $(A, \vdash)$, where $A$ is a $\Sigma$-structure and $\vdash$ is a Tarskian consequence relation on $A$.

A logic is a pair $(\Sigma,\vdash)$ where $\Sigma$ is a signature and $\vdash$ is a Tarskian consequence relation on $\Fm{\Sigma}{X}$.
\end{Df}

A consequence determines, and is determined by, a closure operator on the power set $\wp(A)$, and we will switch freely between the two.

\begin{Df}[Filters, matrices] \label{filter-df} Let $l=(\Sigma,\vdash)$ be a logic.
\begin{enumerate}
\item Let  $A\in \Sigma\text{-Str}$. A subset $F \subseteq A$ is an $l$-filter if for every $\Gamma\cup\{\varphi\}\subseteq \Fm{\Sigma}{X}$ such that $\Gamma\vdash \varphi$ and every valuation $v:\Fm{}{X}\to A$, if $v[\Gamma]\subseteq F$ then $v(\varphi)\in F$. 

In particular, since $l$ is structural,  when $A=\Fm{}{X}$, then the $l$-filters in $\Fm{}{X}$ coincide with the theories of $l$, i.e., the  $\vdash$-closed subsets of $\Fm{}{X}$.

\item Let  $A\in \Sigma\text{-Str}$ and $F \subseteq A$ be an $l$-filter. Then the pair $\langle M,F\rangle$ is then said to be a matrix model of $l$, or shortly, an $l$-matrix. 
\end{enumerate}
\end{Df}

\begin{Obs}[Functoriality of filters]\label{RemarkFunctorialityOfFilters}
Let $l=(\Sigma,\vdash)$ be a logic and $A$ be a $\Sigma$-structure.  Denote by $(Fi_{l}(A), \subseteq) \overset{i_A}\hookrightarrow (P(A), \subseteq)$ the (sub)poset of all $l$-filters of $A$. Clearly $Fi_l(A)$ is closed under arbitrary intersections and, since $l$ is a finitary logic, $Fi_l(A)$ is closed under directed unions. Thus $(Fi_{l}(A), \subseteq)$ is an algebraic lattice where the compact elements are the finitely generated filters. Further, for a homomorphism of $\Sigma$-algebras $f \colon A \to A'$
 and a filter $F' \in Fi_{l}(A')$, the preimage $f^{-1}(F')$ is a filter in $Fi_{l}(A)$. Defining $Fi_{l}(f) := f^{-1}$ (inverse image), then $Fi_{l}$ becomes a contravariant functor from the category $\Sigma\text{-Str}$ to the category $AL$.
\end{Obs}

\subsection{Congruences and equational consequence}\label{SubsectionCongruencesAndEquationalConsequence}

A class ${\bf K} \subseteq \Sigma\text{-Str}$ is a quasivariety if it  is axiomatizable by quasi-identities in the signature $\Sigma$, i.e., formulas of the form
\[(p_{1}\approx q_{1}\wedge...\wedge p_{n}\approx q_{n})\to p\approx q.\]
Quasivarieties are precisely the classes of $\Sigma$-structures that are closed under isomorphisms, subalgebras, products and directed colimits.

A ${\bf K}$-congruence on a $\Sigma$-structure is a $\Sigma$-congruence $\theta \subseteq A \times A$ such that $A/\theta\in \bf K$.


\begin{Obs}[functoriality of ${\bf K}$-congruences]
Let ${\bf K}$ be a quasivariety. For a $\Sigma$-structure $A$ denote by $(Co_{\bf K}(A), \subseteq) \overset{i_A}\hookrightarrow (P(A \times A), \subseteq)$ the (sub)poset of all ${\bf K}$-congruences. As ${\bf K}$ is closed under isomorphisms, subalgebras and products, $Co_{\bf K}(A) \subseteq P(A \times A)$ is closed under arbitrary intersections and, since ${\bf K}$ is closed under directed colimits, $Co_{\bf K}(A) \subseteq P(A \times A)$ is closed under directed unions. Thus $(Co_{\bf K}(A), \subseteq)$ is an algebraic lattice where the compact elements are its finitely generated members.
If $f\colon A \to A'$ is a homomorphism and $\theta' \in Co_{\bf K}(A')$ then, since ${\bf K}$ is closed under isomorphisms and subalgebras, $(f\times f)^{-1}(\theta') \in Co_{\bf K}(A)$. Defining $Co_{\bf K}(f) := (f \times f)^{-1}$, then $Co_{\bf K}$ becomes a contravariant functor from the category $\Sigma\text{-Str}$ to the category $AL$.
\end{Obs}

Let $\Sigma$ be a signature and ${\bf K} \subseteq \Sigma\text{-Str}$ a class of algebras. Let $Eq:=\Fm{}{X}_{\{x\}}\times \Fm{}{X}_{\{x\}}$ be the set of pairs of formulas in at most the variable $x$. We think of such a pair of formulas $\langle \delta, \epsilon \rangle$ as an equation between the terms of both sides and write, following Blok-Pigozzi, $\delta(x)\approx\epsilon(x)$.

\begin{Df}[equational consequence]\label{DefEquationalConsequence}
Given a class of algebras $\textbf{K}$ over the signature $\Sigma$, and a set of variables $X$, the associated equational consequence is the relation $\models_{\textbf{K}}^X$ between a set of equations $\Gamma$ and a single equation $\varphi\approx \psi$  over $\Sigma$ defined by:
\[\Gamma\models_{\textbf{K}}^X\varphi\approx \psi\ i\!f\!\!f\ for\ every\ A\in\textbf{K}\ and\ every\ \Sigma-homomorphism \ \ h:\Fm{}{X}\to A,\]
\[\ if \ h(\eta)=h(\nu)\ for\ all\ \eta\approx\nu\in \Gamma,\ then\ h(\varphi)=h(\psi).\]
\end{Df}
%
In the context of amalgamation we will consider equational consequence for several sets of variables simultaneously. The following is the basic relation that we will need.

\begin{Lem}\label{LemmaInclusionOfVariablesGivesConservativeTranslationOfEquationalConsequences}
Let $\textbf{K}$ be a class of algebras over the signature $\Sigma$, $X \subseteq Y$ sets of variables, and $\Gamma \cup\{\varphi\approx \psi\}$ sets of equations with variables in $X$. Then 
\[\Gamma\models_{\textbf{K}}^Y\varphi\approx \psi\ \ \ \Leftrightarrow\ \ \ \Gamma\models_{\textbf{K}}^X\varphi\approx \psi\]
\end{Lem}
\Dem
$\Leftarrow \colon$ Suppose $\Gamma\models_{\textbf{K}}^X\varphi\approx \psi$.
Let $A\in\textbf{K}$ and $h\colon \Fm{}{Y}\to A$ be a $\Sigma$-homomorphism such that $h(\eta)=h(\nu)$ for all $\eta\approx\nu\in \Gamma$. We need to show that $h(\varphi)=h(\psi)$. For this, just precompose with the inclusion $i \colon \Fm{}{X}\hookrightarrow \Fm{}{Y}$. As $\Gamma \cup\{\varphi\approx \psi\} \subseteq \Fm{}{X}$, we have $(h \circ i)(\eta)=(h \circ i)(\nu)$ for all $\eta\approx\nu\in \Gamma$. From the assumption $\Gamma\models_{\textbf{K}}^X\varphi\approx \psi$ it follows that $h(\varphi)=(h \circ i)(\varphi)=(h \circ i)(\psi)=h(\psi)$.

$\Rightarrow \colon$ Suppose $\Gamma\models_{\textbf{K}}^Y\varphi\approx \psi$.
Let $A\in\textbf{K}$ and $h\colon \Fm{}{X}\to A$ be a $\Sigma$-homomorphism such that $h(\eta)=h(\nu)$ for all $\eta\approx\nu\in \Gamma$. We need to show that $h(\varphi)=h(\psi)$.

Choose a left inverse $Y \to X$ to the inclusion $X \subseteq Y$. These two maps induce homomorphisms  $i \colon \Fm{}{X}\hookrightarrow \Fm{}{Y}$ and $r \colon \Fm{}{Y}\to \Fm{}{X}$ such that $r \circ i = id_{\Fm{}{X}}$. The composition $h \circ r \colon \Fm{}{Y}\to\Fm{}{X}\to A$ is a homomorphism such that $(h \circ r)(\eta)=(h \circ r)(\nu)$ for all $\eta\approx\nu\in \Gamma$. As $\Gamma\models_{\textbf{K}}^Y\varphi\approx \psi$, we have that $(h \circ r)(\varphi)=(h \circ r)(\psi)$. Precomposing with $i$ yields $h(\varphi)=(h \circ r \circ i)(\varphi)=(h \circ r \circ i)(\psi)=h(\psi)$.
\qed

It is shown in chapter 2 of \cite{BP1} that $\vDash_{\bf K}^X$ satisfies the axioms of a structural consequence relation, which is finitary if and only if $\vDash_{\bf K}^X = \vDash_{{\bf K}^Q}^X$, where ${\bf K}^Q$ denotes the quasivariety generated by ${\bf K}$. As usual, the consequence relation $\vDash_{\bf K}^X$ corresponds to a closure operator $Cn_{\bf K}^X$ on $\wp(Eq)$.  A set of equations $\Gamma \subseteq Eq$ is called a ${\bf K}$-theory if is closed under ${\bf K}$-consequence, i.e. if $Cn_{\bf K}^X(\Gamma) = \Gamma$.

For the following Lemma keep in mind that for a set $\Gamma \subseteq Eq$ we have two view points: we can see $\Gamma$ as a set of equations, or as a set of pairs of formulas. From the latter view point it makes sense to ask whether a set is a congruence relation or not. We suppose all the formulas to be in $\Fm{}{X}$ for some set of variables $X$.

\begin{Lem}\label{LemmaEquationalClosureEqualsGeneratedCongruence}
 For a set of equations $\Gamma=\{\langle \gamma^i_1, \gamma^i_2 \rangle \mid i \in I \} \subseteq Eq$ its closure $Cn_{\bf K}^X(\Gamma)$ is the congruence relation generated by $\Gamma$.
 In particular a set $\Gamma \subseteq Eq$ is a ${\bf K}$-theory if and only if is a congruence relation.
\end{Lem}
\Dem
Denote the ${\bf K}$-consequence relation generated by $\Gamma$ by $\hat{\Gamma}$.

$Cn_{\bf K}^X(\Gamma) \supseteq \hat{\Gamma}$: Let $\langle \varphi(x), \psi(x) \rangle \in \hat{\Gamma}$. Note that $\hat{\Gamma}$ is the smallest congruence $\theta$ on $\Fm{}{X}$ such that $\bar{\gamma^i_1}=\bar{ \gamma^i_2}$ in $\Fm{}{X}/\theta$. Thus, every quotient $\Fm{}{X}/\theta$ of $\Fm{}{X}$ in which $\bar{\gamma^i_1(x)}=\bar{ \gamma^i_2(x)}\ (i \in I)$ holds will be itself a quotient of $\Fm{}{X}/\hat{\Gamma}$.

Now suppose that for some algebra $A \in {\bf K}$ and some $a\in A$ we have $(\gamma^i_1)^A(a)= (\gamma^i_2)^A(a)$. Then there is a homomorphism $h \colon \Fm{}{X} \to A$ sending $x$ to $a$ and the other variables to arbitrary elements. The image $Im\, h$ of this homomorphism is isomorphic to a quotient of $\Fm{}{X}$ and a subalgebra of $A$ in which $(\gamma^i_1)^A(a)= (\gamma^i_2)^A(a)$ holds. Thus $Im\, h$ is also isomorphic to a quotient of $\Fm{}{X}/\hat{\Gamma}$, and since $\langle \varphi(x), \psi(x) \rangle \in \hat{\Gamma}$, we have $\varphi^A(a)= \psi^A(a)$ in $Im\, h \subseteq A$. This shows $\langle \varphi(x), \psi(x) \rangle \in Cn_{\bf K}(\Gamma)$.

$Cn_{\bf K}^X(\Gamma) \subseteq \hat{\Gamma}$: Let $\langle \varphi(x), \psi(x) \rangle \in Cn_{\bf K}^X(\Gamma)$. Then since in $\Fm{}{X}/\hat{\Gamma}$ we have $\bar{\gamma^i_1(x)}=\bar{ \gamma^i_2(x)}\ (i \in I)$, we also have $\bar{\varphi(x)} = \bar{\psi(x)}$. Therefore $\langle \varphi(x), \psi(x) \rangle \in \hat{\Gamma}$.
\qed

\subsection{Algebraicity of logics}




We summarize several ways in which a logic can be tied to, and its properties determined by, algebra.

\begin{Df}[algebraic semantics]
 A logic $l=(\Sigma, \vdash_l)$ has an algebraic semantics in a class of algebras ${\bf K}$ given by a set of equations $\delta(x) \approx \epsilon(x)$ if the following equivalence holds: $$\Gamma \vdash_l \varphi \ \ \ \ \Leftrightarrow \ \ \ \ \{\delta(\gamma) \approx \epsilon(\gamma) \mid \gamma \in \Gamma\} \vDash_{{\bf K}} \delta(\varphi) \approx \epsilon(\varphi)$$
\end{Df}

In the following ``algebraizable logic'' will mean ``algebraizable''  in the sense of Blok-Pigozzi \cite{BP1}.

\begin{Df}[algebraizable] \label{BP-def}
A logic $l=(\Sigma, \vdash_l)$ is algebraizable with equivalent algebraic semantics $\textbf{K}$ if
\begin{enumerate}
\item[(1)] $l$ has an algebraic semantics in ${\bf K}$ given by a set of equations $\delta(x) \approx \epsilon(x)$ 
\item[(2)]there is a finite set $\Delta_{j}(p,q),j=1,...,m$ of formulas in two variables such that for every set of equations $\Gamma$ and for every equation $\varphi\approx\psi$, we have that \[\Gamma\models_{\bf K}\varphi\approx\psi\ \Leftrightarrow\ \{\xi\Delta\eta:\xi\approx\eta\in\Gamma\}\vdash_l\varphi\Delta\psi\]
\item[(3)]For every $\psi\in \Fm{}{X}$ we have that
\[\psi\dashv  \vdash_l\Delta(\tau(\psi)).\]
\end{enumerate}
\end{Df}

%
%
%
%
%
%
%

\begin{Df}[compatible congruence, Leibniz operator, reduced matrix]\label{DefCompatibleCongruence} Let $\Sigma$ be a signature, ${A}$ be a $\Sigma$-algebra  and $F\subseteq A$.

(a) Let $\theta$ be a congruence in ${A}$. $\theta$ is said to be compatible with $F$ if, for all $a,b\in A$, if $a\in F$ and $\langle a,b\rangle\in\theta$  then   $b\in F$. Given an algebra $A$ and a subset $F$ of its domain there always exists a greatest congruence of $A$ compatible with $F$ (Theorem 1.5, \cite{BP1}), denoted by $\Omega^{A}(F)$. The function $\Omega^{A}$ with domain the set of all subsets of $A$ is called the Leibniz operator on $A$.


(b) A matrix $(A,F)$ for a logic $l = (\Sigma, \vdash)$ is reduced when $\Omega^{A}(F) = Diag(A) = \{(a,b) \in A \times A : a=b\}$. We will denote $Matr^*_l$ the class of all reduced matrices for $l$; $Alg^{*}_l$ stands for the class of all $\Sigma$-algebras underlying to some reduced matrix.
\end{Df}


\begin{Teo}(The Isomorphism Theorem, first version \cite[Thm 5.1]{BP1}).\label{IT2} \label{TheoremFirstIsoThmInclUniquenessOfTheQuasivariety}

Let $l$ be a logic and ${\bf K}$ a quasivariety. The following are equivalent.
\begin{enumerate}
\item $l$ is algebraizable with equivalent semantics ${\bf K}$
\item For every algebra $A$ the Leibniz operator $\Omega^{A} : Fi_{l}(A) \to Co_{{\bf K}}(A)$ is an isomorphism of algebraic lattices.
\end{enumerate}

Moreover, the quasivariety ${\bf K}$ is unique for providing an equivalent semantics by \cite[Thm 2.15]{BP1}.
\end{Teo}

\begin{Teo}(The Isomorphism Theorem, 2nd version \cite[Thm. 3.58]{Fon}).\label{IT1} \label{TheoremSecondIsoThmInclUniquenessOfTheQuasivariety}

Let $l$ be a logic and ${\bf K}$ be a quasivariety. The following conditions are equivalent:
\begin{enumerate}[{(i)}]
\item $l$ is algebraizable with equivalent algebraic semantics the class ${\bf K}$.
\item For every algebra $A$ there is an isomorphism $\Phi^{A}$ between the algebraic lattices $Fi_{l}(A)$ and $Co_{{\bf K}}(A)$ that commutes
with endomorphisms, i.e., for every $F\in Fi_{l}(A)$ and every $h\in End(A)$, $\Phi^{A}h^{-1}(F)=h^{-1}\Phi^{A}F$.
\item There is an isomorphism of algebraic lattices $\Phi^{Fm} \colon \mathcal{T}h(l) \to Co_{{\bf K}}(Fm)$ that commutes with
substitutions, i.e., for every $T\in \mathcal{T}h(l)$ and every $\sigma\in End(Fm)$, $\Phi\sigma^{-1}T=\sigma^{-1}\Phi T$.
\end{enumerate}
If the above conditions are satisfied, the isomorphism of (iii) is unique and given by the Leibniz operator $\Omega$. It can be expressed by via set $\Delta(p,q)$ of formulas in two variables as $\Omega(F)=\{\langle \varphi, \psi \rangle \mid \Delta(\varphi, \psi) \subseteq F \}$ and its inverse can be expressed via a set of equations $\delta(x)=\epsilon(x)$ as $\Omega^{-1}(\theta)=\{ \varphi \in Fm \mid \langle \delta(x),\epsilon(x)\rangle \subseteq \theta \}$.

For a more general algebra $A$, the same set of formulas, resp. set of equations, induce isomorphisms and there is at most one way to express an isomorphism by equations in this way (\cite[Exercise 3.40]{Fon}).
\end{Teo}


The following corollary was the main motivation to define the congruence filter pairs in order to establish conditions to get algebraizable logics knowing previously a quasivariety of some similarity type $\Sigma$. 

\begin{Cor}(The Isomorphism Theorem, 3rd version)

Let $l$ be a logic and ${\bf K}$ be a quasivariety. Then following conditions are equivalent:
\begin{enumerate}
\item $l$ is an algebraizable logic with equivalent algebraic semantics the class ${\bf K}$.
\item There is a natural isomorphism between the functors $Fi_{l}$ and $Co_{{\bf K}}$.
\end{enumerate}
\end{Cor}

\Dem

$''1\Rightarrow 2''$ Suppose that $l$ is an algebraizable logic. By theorem \ref{IT2} we have that for every $A\in\Sigma\text{-Str}$, the Leibniz operator $\Omega^{A}:Fi_{l}(A)\to Co_{K}(A)$ is a isomorphism. Let $\Omega=(\Omega^{A})_{A\in\Sigma\text{-Str}}$. We prove that $\Omega$ is a natural transformation. In other to do that, it is enough to prove that given $f\in Hom_{\Sigma\text{-Str}}(A,B)$, the following diagram commutes

\[
\xymatrix{
Fi_{l}(A)\ar[r]^{\Omega^{A}}&Co_{{\bf K}}(A)\\
Fi_{l}(B)\ar[u]^{f^{-1}}\ar[r]_{\Omega^{B}}&Co_{{\bf K}}(B)\ar[u]_{{(f \times f)}^{-1}}
}\]

Let $F\in Fi_{l}(B)$. Firstly we prove that $f^{-1}(\Omega^{B}(F))$ is compatible with $f^{-1}(F)$. Let $(a,b)\in f^{-1}(\Omega^{B}(F))$ and suppose that $a\in f^{-1}(F)$.  Then $(f(a),f(b))\in \Omega^{B}(F)$ and $f(a)\in F$. Therefore $f(b)\in F$, thus $b\in f^{-1}(F)$. Thus $(a,b)\in \Omega^{A}((f \times f)^{-1}(F))$.

Now let $(a,b)\in \Omega^{A}(f^{-1}(F))$, then by algebraizability of $l$, we have $\Delta^{A}(a,b)\in f^{-1}(F)$. Thus \linebreak $\Delta^{B}(f(a),f(b)) = f(\Delta^{A}(a,b))\subseteq F$. Therefore $(f(a),f(b))\in \Omega^{B}(F)$ and finally $(a,b)\in f^{-1}(\Omega^{B}(F))$. So $\Omega^{A}(f^{-1}(F))=f^{-1}(\Omega^{B}(F))$, which means that $\Omega$ is natural.
\vtres

$''2\Rightarrow 1''$ Suppose that there is a natural isomorphism $\Phi:Fi_{l}\to Co_{{\bf K}}$. In particular we have that for every $A\in\Sigma\text{-Str}$, $\Phi^{A}:Fi_{l}(A)\to Co_{{\bf K}}(A)$ is a isomorphism and commutes with endomorphisms. By theorem \ref{IT1} we have that $l$ is an algebraizable logic.
\qed
\vtres

We consider now some special kinds of logics that are generalizations for algebraizable logics (see \cite{Fon}).

\begin{Df}[protoalgebraic, truth-equational]\label{SpecialLogics} 
Let $l=(\Sigma,\vdash)$ a logic:
\begin{itemize}
\item $l$ is a {protoalgebraic} logic if for any theory $T\in Th(l)$,
\[ if \langle\varphi,\psi\rangle\in\Omega(T)\ then\ T,\varphi\vdash\psi\ and\ T,\psi\vdash\varphi.\]
\item A class of matrices $M$ has its filters equationally definable by a set of equations $\tau(p)$ if for
every matrix $\langle A,F\rangle\in M$, for every $a\in A$,
\[a\in F\ \ iff\ \ \delta^{A}(a) = \varepsilon^{A}(a),\ for\ every\ \delta\approx\varepsilon\in \tau(p).\]
$l$ is a {truth-equational} logic if the class of reduced matrix $Matr^{*}_l$ has its filters equationally definable.

\item Let $K$ be a pointed class of algebra of a similarity type $\Sigma$. $l$ is the assertional logic of $K$ if,
\[\Gamma\vdash\varphi \Leftrightarrow \{\gamma\approx T\ :\ \gamma\in\Gamma\}\models_{K}\varphi\approx
 T\]
\end{itemize}
\end{Df}

Now we recall some characterizations of the classes of logics defined above.

\begin{Teo}\label{teocar}
Let $l$ be a logic:
\begin{itemize}
\item $l$ is protoalgebraizable iff $\Omega$ is monotone on set of theories $Th(l)$.
\item $l$ is equivalential iff $(\Omega^{A})_{A\in\Sigma\text{-Str}}$ commutes with homomorphism (categorial naturality) and $\Omega$ is monotone.
\item $l$ is truth-equational iff there exists a set of equations $\tau(p)$ such that for every algebra $A$ and every $F\in Fi_{l}(A)$,
\[F=\{a\in A;\ \tau^{A}(a)\subseteq\Omega^{A}(F)\}\]
\end{itemize}
\end{Teo}

The following diagram summarizes the relations between the different classes of logics mentioned in this section. 

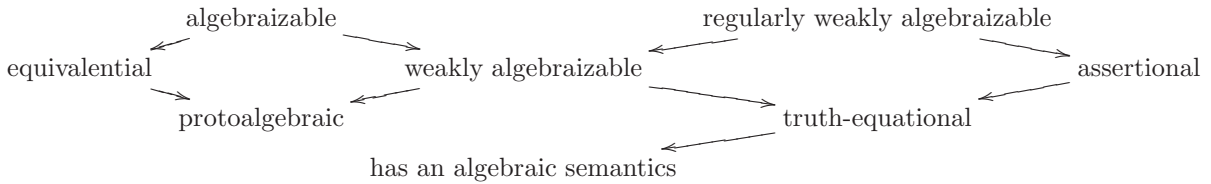
\[
\xymatrix@=0.4em{
&\text{algebraizable}\ar[dr]\ar[dl]&& \text{regularly weakly algebraizable} \ar[dl]\ar[rd]&\\
\text{equivalential}\ar[dr]&& \text{weakly algebraizable}\ar[dl]\ar[dr]&&\text{assertional}\ar[ld]\\
&\text{protoalgebraic} &&\text{truth-equational} \ar[dl]\\
&& \text{has an algebraic semantics} &
}\]
%








\section{Finitary Filter Pairs}\label{SectionFilterFunctors}

To motivate the notion of filter pair, we recall from Remark \ref{RemarkFunctorialityOfFilters} that, given a logic $l=(\Sigma, \vdash)$, one can associate to each $\Sigma$-structure $M$ its collection of $l$-filters $Fi_{l}(M)$ and that, together with taking inverse image, this constitutes a functor $Fi_l \colon \Sigma\text{-Str} \to AL$, from $\Sigma$-structures to algebraic lattices. Also recall from Remark \ref{RemarkFunctorialityOfFilters} that for $M\in \Sigma\text{-Str}$ the inclusion $i_{M} :  (Fi_{l}(M), \subseteq) \hookrightarrow (\mathcal{P}(M),\subseteq)$ preserves arbitrary infima and directed suprema, i.e., it is a morphism in the category $AL$ (in particular preserves order).

Moreover, given a morphism $h:M\to N$ we have the following commutative diagram:

\[
\xymatrix{
M\ar[d]_{h}&Fi_{l}(M)\ar@{^(->}[r]^{i_{M}}&(\mathcal{P}(M),\subseteq)\\
N&Fi_{l}(N)\ar[u]^{h^{-1}}\ar@{^(->}[r]_{i_{N}}&(\mathcal{P}(N), \subseteq)\ar[u]_{h^{-1}}\\
}
\]

This collection of data is the motivating example for the notion, and the name, of filter pair. 

\begin{Df}\label{DefinitionFilterPair}
Let $\Sigma$ be a signature. A {\bf filter pair} over $\Sigma$ is a pair $(G,i)$, consisting of a contravariant functor $G:\Sigma\text{-Str}\to AL$ and a collection of maps $i=(i_{M})_{M\in \Sigma\text{-Str}}$ such that for any $M\in\Sigma\text{-Str}$ the function $i_{M}:G(M)\to (\mathcal{P}(M),\subseteq)$ satisfies  the following properties:

{\bf 1.} For any $M\in\Sigma\text{-Str}$, $i_{M}$ preserves arbitrary infima (in particular $i_{M}(\top)=M$) and directed suprema.

{\bf 2.} Given a morphism $h:M\to N$ the following diagram commutes:

\[
\xymatrix{
M\ar[d]_{h}&G(M)\ar[r]^{i^{G}_{M}}&(\mathcal{P}(M);\subseteq)\\
N&G(N)\ar[u]^{G(h)}\ar[r]_{i^{G}_{N}}&(\mathcal{P}(N);\subseteq)\ar[u]_{h^{-1}}
}
\]
\end{Df}

\begin{Obs}
 Condition {\bf 2.} says that $i$ is a natural transformation from $G$ to the functor $\wp\colon \Sigma\text{-Str}^{op}\to AL$ sending a $\Sigma$-structure to the power set of its underlying set and a homomorphism of $\Sigma$-structures to its associated inverse image function.
\end{Obs}

From a filter pair we obtain a consequence relation on every $\Sigma$-structure, i.e. an abstract logic in the sense of Definition \ref{DefAbstractLogic}.

\begin{Prop}\label{LogicsFromFilterPairs}
Let $(G,i)$ be a filter pair over the signature $\Sigma$ and $A$ a $\Sigma$-structure. Then there is a generalized logic $l_{G}^A=(A,\vdash_{G}^A)$, defined as follows:\\
Given $\Gamma\cup\{\varphi\}\subseteq A$, define
\[\Gamma\vdash_{G}\varphi\ \ i\!f\!\!f\ \ for\ any\ a\in G(\Fm{}{X}),\ if\ \Gamma\subseteq i_{\Fm{}{X}}(a)\ then\ \varphi\in i_{\Fm{}{X}}(a).\]
\end{Prop}

\Dem

It is easy to see that $\vdash_{G}^A$ satisfies reflexivity, cut and monotonicity.

The structurality is a consequence of condition {\bf 2} (naturality). Indeed, let $\sigma\in hom(\Fm{}{X},\Fm{}{X})$ and $\Gamma\cup\{\varphi\}\subseteq \Fm{}{X}$ such that $\Gamma\vdash_{G}\varphi$. Consider $a\in G(\Fm{}{X})$ such that $\sigma[\Gamma]\subseteq i^{G}_{\Fm{}{X}}(a)$. This implies $\Gamma\subseteq\sigma^{-1}(i^{G}_{\Fm{}{X}}(a))$. By naturality we have $\sigma^{-1}(i^{G}_{\Fm{}{X}}(a))=i^{G}_{\Fm{}{X}}(G(\sigma)(a))$. Therefore $\varphi\in i^{G}_{\Fm{}{X}}(G(\sigma)(a))=\sigma^{-1}(i^{G}_{\Fm{}{X}}(a))$ and finally $\sigma(\varphi)\in i^{G}_{\Fm{}{X}}(a)$.

Now we are going to prove the finitarity. Let $\Gamma\cup\{\varphi\}\subseteq \Fm{}{X}$. Consider the set $S=\{\Gamma'\subseteq \Fm{}{X};\ \Gamma'\subseteq_{fin}\Gamma\}$. Notice that $S$ is a directed set. Suppose that for any $\Gamma'\in S,\ \Gamma'\nvdash_{G}\varphi$, hence there is $a\in G(\Fm{}{X})$ such that $\Gamma'\subseteq i_{\Fm{}{X}}(a)$ and $\varphi\not\in i_{\Fm{}{X}}(a)$. Denote by $a_{\Gamma'}=\wedge\{a\in G(\Fm{}{X});\ \Gamma'\subseteq i_{\Fm{}{X}}(a)\}$ . $i_{\Fm{}{X}}$ preserves inf, thus $\Gamma'\subseteq i_{\Fm{}{X}}(a_{\Gamma'})$ and $\varphi\not\in i_{\Fm{}{X}}(a_{\Gamma'})$. We obtain that the set $s=\{a_{\Gamma'};\ \Gamma'\in S\}$ is a directed set.

By the assumption $i_{\Fm{}{X}}$ preserves directed suprema, hence \[\Gamma=\cup S\subseteq \bigcup_{\Gamma'\in S}i_{\Fm{}{X}}(a_{\Gamma'})=i_{\Fm{}{X}}(\vee s).\] On the other hand $\varphi\not\in\bigcup_{\Gamma'\in S'}i_{\Fm{}{X}}(a_{\Gamma'})=i_{\Fm{}{X}}(\vee s)$. Therefore $\Gamma\nvdash_{G}\varphi.$\qed

\begin{Obs}
 The main instance of interest of the construction of Prop. \ref{LogicsFromFilterPairs} is the case $A=\Fm{\Sigma}{X}$, the formula algebra over some set $X$ of variables. In this case one obtains a logic in the usual sense of the word, which we will denote by $(\Sigma, \vdash^X)$. Thus one can see a filter pair (together with a set $X$ of variables) as a \emph{presentation} of a logic, different in style from the usual presentations by axioms and rules or by matrices. It is clear from the definition that the logic defined in this way does not depend on the values of the filter pair at $\Sigma$-structures other than $\Fm{\Sigma}{X}$, and indeed it only depends on \emph{the image of the map} $G(\Fm{\Sigma}{X}) \to \wp(\Fm{\Sigma}{X})$, as this is exactly the collection of theories of the logic, by definition.
 
 Thus, just as an algebraic structure can have many different presentations by generators and relations, a logic can have presentations by different filter pairs, each of which can be useful for different purposes.
\end{Obs}

\begin{Ex}\label{ExampleTheFilterPairOfFilters}
Given a Tarskian logic $l=(\Sigma,\vdash)$, by Remark \ref{RemarkFunctorialityOfFilters}, defining $Fi_{l}(A)$ to be the set of $l$-filters on a $\Sigma$-structure $A$ provides a functor $Fi_{l} : \Sigma\text{-Str}^{op}\to AL$, and hence a filter pair $(Fi_{l},i)$ where and $i$ is the inclusion of filters into all subsets. As filters on the formula algebra are exactly the theories, this shows that every logic admits a presentation by a filter pair.
\end{Ex}

There is a direct description of the closure operator associated to the consequence relation of Prop. \ref{LogicsFromFilterPairs} through an application of the notion of adjunctions between posets. 
Since the natural transformation $i$ of a filter pair $(G,i)$ preserves infima, by Theorem \ref{TheoremExistenceOfLeftAdjoint} for each algebra $A$, the map $i_A$  has a (unique) left adjoint.

\begin{Df}\label{DefLeftAdjointOfAFilterPair}
 Let $(G,i)$ be a filter pair. The left adjoint of $i_A$ will be denoted by $\Xi(i)_A$.
\end{Df}




Note that in general there is no reason that the collection of left adjoints $\Xi(i)_A\ \ (A \in \Sigma\text{-Str})$ should form a natural transformation. It is Nevertheless a resource that comes for free with a filter pair and can be usefully employed to understand the properties of the associated logics.

\begin{Prop}\label{PropDescriptionOfConsequenceRelationByClosureOperator}
 Let $(G,i^{G})$ be a filter pair and $A$ a $\Sigma$-structure. Then the closure operator on $A$ associated to the consequence relation $\vdash_{G}^A$ on $A$ of Prop. \ref{LogicsFromFilterPairs} is the one given by $i_A \circ \Xi(i)_A$.
\end{Prop}

\Dem
Recall the definition of $\vdash_{G}^A$:
\[D \vdash_{G}^A a\ \ i\!f\!\!f\ \ for\ any\ z\in G(A),\ if\ D\subseteq i_{A}(z)\ then\ a\in i_{A}(z).\]

Consider the set $Z_D:=\{z\in G(A) \mid  D\subseteq i_{A}(z) \}$. We can rephrase the definition of $\vdash_G^A$ by saying 
\[D \vdash_{G}^A a\ \ i\!f\!\!f\ \ \forall z\in Z_D,\ a\in i_{A}(z)\]

Since $G(A)$ is complete, there exists an infimum, say $z'$, of $Z_D$. Since $i_A$ preserves arbitrary infima, we have $i_A(z')=\bigcap_{z \in Z_D} i_A(z)$. Since we have $D \subseteq i_A(z) \forall z \in Z_D$, we also have $D \subseteq \bigcap_{z \in Z_D} i_A(z) = i_A(z')$. Therefore $z' \in Z_D$, i.e. it is the minimal $z'$ for which $i_A(z)$ contains $D$. Since $i_A$ is order preserving, if $a \in i_A(z')$, then $a \in i_A(z)$ for every $z \in Z_D$ and hence $D \vdash_{G}^A a$. The opposite direction is also true, i.e. if $D \vdash_{G}^A a$, then $a \in i_A(z')$, just because $z' \in Z_D$.

Thus we can describe the consequence relation by $D \vdash_{G}^A a \ \text{ iff } \ a \in i_{A}(z')$.

By the formula for left adjoints from Theorem \ref{TheoremExistenceOfLeftAdjoint}, for a subset $D \subseteq A$ we have $\Xi(i)_A(D)=\bigwedge \{z \in G(A) \mid D \leq i_A(z) \} = z'$. Altogether we obtain
\[D \vdash_{G}^A a\ \ \text{ iff } \ a \in i_{A}(\Xi(i)_A(D))\]
This means exactly that the closure operator $i_A \circ \Xi(i)_A$ is the closure operator associated to the consequence relation $\vdash_G^A$.
\qed

In the remainder of the section we will explore the relationship between filters of the logic $(\Sigma, \vdash^X)$ associated to a filter pair $(G,i)$, and sets in the image of $i_A$ for a $\Sigma$-structure $A$. We start by giving a name to the latter.

\begin{Df}\label{DefiFilter}
 Let $(G,i)$ be a filter pair and $A$ a $\Sigma$-structure. An $i$-filter in $A$ is a subset  $F \subseteq A$ such that for all $D \subseteq A$ the following implication holds:
\[ ( D\subseteq F\text{ and }D \vdash_{G}^A a)\ \Rightarrow a\in F \]
In other words, $i$-filters in $A$ are the theories for the abstract logic $(A,\vdash)$. 
\end{Df}

\begin{Lem}\label{Lemi-FiltersAreTheImageOfi}
Let $(G,i^{G})$ be a filter pair and $A$ a $\Sigma$-structure. The image of $i_A$ in $\wp(A)$ consists exactly of the $i$-filters in $A$.
\end{Lem}
\Dem
The set $i-Fi(A)$ of Filters in a $\Sigma$-algebra $A$ is defined by the following condition:
\[F \in i-Fi(A) \ \ \text{ iff } \  [ ( D\subseteq F\text{ and }D \vdash_{G}^A a)\ \Rightarrow a\in F]\] 
In view of the previous proposition this means 
\[F \in i-Fi(A) \ \ \text{ iff } \  [  (D\subseteq F\text{ and }a \in (i_A \circ \Xi(i)_A)(D)\ \Rightarrow a\in F]\] 
It is enough to check the condition on the right hand side on the maximal $D \subseteq F$, i.e. on $F$ itself, which leads to 
\[F \in i-Fi(A) \ \ \text{ iff } \  [ a \in (i_A \circ \Xi(i)_A)(F) \Rightarrow a\in F] \]
Now since $i_A \circ \Xi(i)_A$ is a closure operator, this means exactly that $F$ is closed for this operator, i.e. $F=(i_A \circ \Xi(i)_A)(F)$, hence $F$ is in the image of $i_A$. Vice versa, it is clear from the properties of Galois adjunctions that sets in the image of $i_A$ are closed. 
\qed

The notions of $i$-filter for a filter pair and filter for a logic associated to that filter pair are, in general, different. The content of the next proposition is that one notion subsumes the other. It can be shown in relevant examples that this inclusion is strict.

\begin{Prop}\label{PropImageOfiConsistsOfFilters}
Let $(G,i)$ be a filter pair, $X$ a set and $l=(\Sigma,\vdash^X)$ the associated logic with set of variables $X$. Then for any algebra $A$ the subsets in the image of $i_A$ are $l$-filters.
\end{Prop}
\Dem
Let $A$ be an algebra, $F=i_A(x) \subseteq A$ a subset in the image of $i$, $\Gamma \cup \{\varphi\} \subseteq \Fm{}{X}$ formulas with $\Gamma \vdash^X \varphi$ and $v \colon \Fm{}{X} \to A$ a homomorphism with $v(\Gamma) \subseteq F$. We need to show that $\varphi \in F$. For this consider the naturality square

\[
\xymatrix{
\Fm{}{X}\ar[d]_{v}&G(\Fm{}{X})\ar[r]^{i^{G}_{\Fm{}{X}}}&(\mathcal{P}(\Fm{}{X});\subseteq)\\
A&G(A)\ar[u]^{G(v)}\ar[r]_{i^{G}_{A}}&(\mathcal{P}(A);\subseteq)\ar[u]_{v^{-1}}
}
\]

We have $\Gamma \subseteq v^{-1}(F) = v^{-1}(i_A(x)) = i_{\Fm{}{X}}(G(v)(x))$, 

So $v^{-1}(F)$, lying in the image of $i_{\Fm{}{X}}$, is a theory containing $\Gamma$. The assumption $\Gamma \vdash^X \varphi$ then implies $\varphi \in v^{-1}(F)$, i.e. $v(\varphi) \in F$.
\qed

\begin{Prop}\label{PropiMatricesDefineTheSameLogicAsAllMatrices}
Given a filter pair $(G,i)$ and a set $X$, the consequence relation $\vdash$ on $\Fm{}{X}$ defined by
 \begin{center}
$\Gamma\vdash\varphi$ iff for any $\Sigma$-structure $M$, for any $a\in G(M)$ and any valuation $v:\Fm{}{X}\to M$, if $v(\Gamma)\subseteq i_{M}(a)$ then $v(\varphi)\in i_{M}(a)$
 \end{center}
coincides with the consequence relation $\vdash^X$.
\end{Prop}
\Dem
Suppose $\Gamma\vdash\varphi$. Then, taking the identity as valuation, one has that if $\Gamma\subseteq i_{\Fm{}{X}}(a)$ then $\varphi\in i_{\Fm{}{X}}(a)$ for all $a \in G(\Fm{}{X})$. Since the sets in the image of $i_{\Fm{}{X}}$ are precisely the theories of the logic $l:=(\Sigma,\vdash^X)$, this shows $\Gamma \vdash^X \varphi$.

On the other hand suppose $\Gamma\vdash^X\varphi$. Then for any $l$-matrix $(M,F)$ and valuation $v:\Fm{}{X}\to M$ one has that  if $v(\Gamma)\subseteq F$ then $v(\varphi)\in F$. Since by Proposition \ref{PropImageOfiConsistsOfFilters} all subsets in the image of $i_M$ are filters, this implies that the defining condition for $\Gamma\vdash\varphi$ is satisfied.
\qed

\begin{Obs}\label{RemarkAlternativeDescriptionsOFLogicsFromFilterPairs}
\begin{enumerate}[(a)]
\item Call a matrix $(A,F)$ an $i$-matrix if the filter $F$ is in the image of $i_A$. Denoting the collection of all $i$-matrices by $i\text{-}Matr$, and the collection of \emph{all} matrices for the logic $l=(\Sigma, \vdash^X)$ by $l\text{-}Matr$, Proposition \ref{PropImageOfiConsistsOfFilters} can be subsumed as stating an inclusion $i\text{-}Matr \subseteq l\text{-}Matr$. Proposition \ref{PropiMatricesDefineTheSameLogicAsAllMatrices} then says that, although this inclusion can be strict, the two classes of matrices always define the same logic.

\item The inclusions of $i$-filters into all filters, established by Proposition \ref{PropImageOfiConsistsOfFilters} are also easily seen to form a natural transformation. Thus we can consider the natural transformation $i^{G}:G\Rightarrow Fi_{l}$. This exhibits $(Fi_l, i)$ as a weakly terminal filter pair among all filter pairs presenting the same logic $l$. The collection of filter pairs presenting a given fixed logic will be studied in a follow-up work.
\end{enumerate}
\end{Obs}

The following two propositions establish that, although the notions of $i$-filter and $l$-filter (where $l=(\Sigma,\vdash^X)$ is the logic with variables $X$) may differ for general algebras, this is not the case  \emph{for an absolutely free algebra} $\Fm{}{Z}$ (with a set of generators $Z$ possibly different from $X$).

\begin{Prop}\label{2.19}\label{PropHomomorphismsInduceTranslationsAndVariableInclusionsInduceConservativeTranslations}
Let $(G,i)$ be a filter pair over the signature $\Sigma$.
\begin{enumerate}[(i)]
\item For any homomorphism of $\Sigma$-structures $f\colon A\to B$ and $\Gamma\cup\{\varphi\}\subseteq A$:
\[\Gamma\vdash^{A} \varphi\ \Rightarrow\ f[\Gamma]\vdash^{B} f(\varphi).\]
\item For any injective map of sets $f:X\rightarrowtail Y$ and $\Gamma\cup\{\varphi\}\subseteq \Fm{}{X}$:
\[\Gamma\vdash^{X}\varphi\ \Leftrightarrow\ f[\Gamma]\vdash^{Y} f(\varphi).\]
\end{enumerate}
\end{Prop}

\Dem $\emph{(i)}$ Suppose $\Gamma\vdash^{A} \varphi$. Let $z\in G(B)$ such that $f[\Gamma]\subseteq i_{B}(z)$. Then $\Gamma\subseteq f^{-1}(i_{B}(z))=i_{A}(G(f)(z))$. Since $\Gamma\vdash^{A}\varphi$, we have that $\varphi\in i_{A}(G(f)(z))$. Therefore $f(\varphi)\in i_{B}(z)$. As $z$ was arbitrary we have $f[\Gamma]\vdash^{B} f(\varphi)$.
\vtres

$\emph{(ii)}$ Let $f:X\to Y$ be injective. By $\emph{1}$ we have that $\Gamma\vdash^{X}\varphi\ \Rightarrow\ f[\Gamma]\vdash^{Y}f(\varphi)$. It remains to prove the converse. Let $z\in G(\Fm{}{X})$ such that $\Gamma\subseteq i^{G}_{X}(z)$. Since $f$ is injective there is a $g\colon Y\to X$ such that $g\circ f=Id_{X}$. Hence $g\circ f[\Gamma]=\Gamma$. Therefore $f[\Gamma]\subseteq g^{-1}(i_{X}(z))=i_{Y}(G(g)(z))$. Since $f[\Gamma]\vdash^{Y}f(\varphi)$, we have $f(\varphi)\in i^{G}_{Y}(G(g)(z))=g^{-1}(i^{G}_{X}(z))$. Therefore $\varphi=g(f(\varphi))\in i^{G}_{X}(z)$.
\qed

\begin{Prop}\label{FreeAlgebraResult}
Let $(G,i)$ be a filter pair and $l:=(\Fm{}{X}, \vdash^X)$ the associated logic with an infinite set of variables $X$. Let $Z$ be another set and $F\subseteq \Fm{}{Z}$ an $l$-filter. Then $F=i_Z(\Xi(F))$. In particular all filters in absolutely free algebras lie in the image of $i$.
\end{Prop}
\Dem
From the general properties of adjunctions we know $F \subseteq i_Z(\Xi(F))$. It remains to show $F \supseteq i_Z(\Xi(F))$. So let $\varphi \in i_Z(\Xi(F))$. By Prop. \ref{PropDescriptionOfConsequenceRelationByClosureOperator} this means $F \vdash^Z \varphi$.

By Prop. \ref{LogicsFromFilterPairs} the logic $(\Fm{}{Z}, \vdash^Z)$ is finitary, so there is a finite subset $F' \subseteq F$ such that $F' \vdash^Z \varphi$. Then there is a finite set $V \subseteq Z$ of variables occurring in the formulas of $F' \cup \{\varphi\}$. Choose a map $h \colon Z \to X$ which is injective on $V$ (this is possible since $X$ is infinite) and a map $v \colon X \to Z$ such that $v \circ h|_V = id_V$. We continue denoting the induced maps on formula algebras by $h$, resp. $v$. By construction we then have $v(h(F))=F$ and $v(h(\varphi))=\varphi$.

From $F' \vdash^Z \varphi$ and Prop. \ref{PropHomomorphismsInduceTranslationsAndVariableInclusionsInduceConservativeTranslations}(i) it follows that $h(F') \vdash^X h(\varphi)$, i.e. $h(F')$ implies $h(\varphi)$ in the logic $l$.

Now $v$ is a valuation such that $v(h(F'))=F'\subseteq F$ and $F$ was supposed to be an $l$-filter, so $\varphi = v(h(\varphi)) \in F$.
\qed

\section{Congruence filter pairs 
}\label{SectionFilterPairsOverCo}

Motivated by the example of algebraizable logics, we introduce in this section a special class of filter pairs, namely the filter pairs whose functor part is given by $Co_{{\bf K}}: \Sigma\text{-Str}\to AL$, where ${\bf K}\subseteq  \Sigma\text{-Str}$ is a class of $\Sigma$-structures, and for $M \in \Sigma\text{-Str}$ $Co_{{\bf K}}(M)$ denotes the ordered set of congruences whose associated quotient lies in ${\bf K}$. We will call such filter pairs \emph{congruence filter pairs}.

We begin by studying the question which classes ${\bf K}\subseteq  \Sigma\text{-Str}$ are eligible for the formation of such a filter pair, i.e. when $Co_{{\bf K}}$ becomes a functor with values in algebraic lattices.







We will provide some  examples of congruence filter pairs and present results relating  these filter pairs with classes of logics in the Leibniz hierarchy. 

\subsection{Functors of congruences}

To begin, we record some conditions on a class ${\bf K} \subseteq \Sigma\text{-Str}$ for obtaining a functor $Co_{\bf K} \colon \Sigma\text{-Str} \to Ord$, where $Ord$ denotes the category of partially ordered sets and order preserving maps. 

\begin{Prop}\label{PropCoKIsAFunctor}
 Let $\Sigma$ be a signature and ${\bf K} \subseteq \Sigma\text{-Str}$ a class of algebras that is closed under subalgebras and isomorphisms. Then there is a functor $Co_{\bf K} \colon \Sigma\text{-Str}^{op} \to Ord$ associating to a $\Sigma$-structure $A$ the set of congruences $\theta$ such that $A/\theta \in {\bf K}$ and to a homomorphism $f \colon A \to B$ the function $$Co_{\bf K}(f)\colon Co_{\bf K}(B) \ni \theta \mapsto (f \times f)^{-1}(\theta)=\{\langle a,a' \rangle \mid \langle f(a),f(a') \rangle \in \theta \} \in Co_{\bf K}(A)$$

 Conversely, if ${\bf K} \subseteq \Sigma\text{-Str}$ is a class of algebras such that $Co_{\bf K} \colon \Sigma\text{-Str}^{op} \to Ord$ is a functor, which on morphisms is given by $Co_{\bf K}(f)=(f \times f)^{-1}$, then ${\bf K}$ is closed under substructures.
\end{Prop}
\Dem
The set-theoretic preimage of a congruence is a congruence again: Indeed, congruences are exactly the set-theoretic preimages of the equality relation (Denoted $\Delta$) along homomorphisms. So, if $\theta \in Co_{\bf K}(B)$ is of the form $(g \times g)^{-1}(\Delta_C)$ for some $g \colon B \to C$, then $(f \times f)^{-1}(\theta)=((g \times g) \circ (f\times f))^{-1}(\Delta_C)$, and hence is a congruence again.

To see that the quotient by this congruence lies in ${\bf K}$, suppose that $\theta \in Co_{\bf K}(B)$. From the isomorphism theorems of universal algebra we obtain a commutative square
\[
\xymatrix{
A \ar@{->>}[d] \ar[r]^{f}&B \ar@{->>}[d]\\
A/f^{-1}(\theta)\ar@{>->}[r]_{\bar{f}}& B/\theta
}
\]
where the lower horizontal map is injective. So $A/f^{-1}(\theta)$ is isomorphic to a subalgebra of $B/\theta$, hence lies in ${\bf K}$ again.

\medskip

Conversely, suppose that $Co_{\bf K}$ is a functor as described and let $B \in {\bf K}$ and $h \colon A \hookrightarrow B$ a substructure. Then $\Delta_B \subseteq B \times B$ is in $Co_{\bf K}(B)$. Hence so is $(h \times h)^{-1}(\Delta_B)=\Delta_A$, hence $A=A/\Delta_A \in {\bf K}$.
\qed

The following result is essentially a result attributed to Malcev (see \cite{Art} Theorem 1.37, pp 869):

\begin{Prop}\label{PropKClosedUnderProducts}
 Let $\Sigma$ be a signature and ${\bf K} \subseteq \Sigma\text{-Str}$ a class of algebras that is closed under subalgebras and isomorphisms and $A \in \Sigma\text{-Str}$. 
 Are equivalent:
 
 (i) $Co_{\bf K}(A)$ has a minimal element;
 
 (ii) ${\bf K}$  is 
 closed under products;
 
 (iii) ${\bf K}$ is a reflective subcategory of $\Sigma\text{-Str}$.

\end{Prop}

\qed

\begin{Cor}\label{CorollaryUnderVopenkasPrincipleKMustBeAQuasivariety}
 Assuming Vop{\v e}nka's principle, under the conditions of Prop. \ref{PropKClosedUnderProducts} the class ${\bf K}$ is a quasivariety.
\end{Cor}
\Dem
Assuming Vop{\v e}nka's principle, by the Corollary after Theorem 6.14 in \cite{AR}, a class of $\Sigma$-structures is a quasivariety if and only if it is closed under subalgebras and products.
\qed

In view of Cor. \ref{CorollaryUnderVopenkasPrincipleKMustBeAQuasivariety} it seems possible that the question whether under the conditions of Prop. \ref{PropKClosedUnderProducts} ${\bf K}$ must be a quasivariety is undecidable. However, for our application to filter pairs we want $Co_{\bf K}$ to be taking values in algebraic lattices and can ask whether this condition forces ${\bf K}$ to be a quasivariety. We have not been able to settle this question and leave it as a question.

\begin{Que}
 Let $\bf K$ be a class of $\Sigma$-structures that is closed under subalgebras and isomorphisms and suppose that the functor $Co_{\bf K}$ of Prop. \ref{PropCoKIsAFunctor} takes values in algebraic lattices. Is $\bf K$ then necessarily a quasivariety?
\end{Que}

\begin{Obs}
 One way of approaching the above question is to note that the proof of Proposition \ref{PropKClosedUnderProducts} shows that, in the situation in question ${\bf K}$ is a regular epi-reflective subcategory of $\Sigma\text{-Str}$, i.e. a reflexive subcategory for which the unit morphisms of the adjunction are regular epimorphisms. It is a consequence of this, as shown in the proof, that ${\bf K}$ is closed under subalgebras and products.
 Quasivarieties are precisely the regular epi-reflective subcategories of $\Sigma\text{-Str}$ which are also closed under directed colimits (see e.g. \cite{AdamekRosickyVitale}, Cor. 10.24 together with the remarks on the next page). A quick proof of the relevant direction proceeds by noting that all we need for having a quasivariety is closure under ultraproducts, and ultraproducts are filtered colimits of pdiagrams of products.
 The question can thus be reformulated as asking whether the additional assumption that $Co_{\bf K}$ takes values in algebraic lattices forces ${\bf K}$ to be closed under directed colimits, and hence to be a quasivariety.
 
 To see the possible connection, note that the reflection functor is given by taking certain quotients. From this we obtain a reflection of ordered sets congruences, i.e. a left inverse $r \colon Co_{\Sigma\text{-Str}}(A) \to Co_{\bf K}(A)$ to the inclusion $Co_{\bf K}(A) \hookrightarrow Co_{\Sigma\text{-Str}}(A)$: given a congruence $\theta \in Co_{\Sigma\text{-Str}}(A)$, we define $r(\theta):=(p \times p)^{-1}(min_{A/\theta}) \in Co_{\bf K}(A)$ where $p$ denotes the projection $p \colon A \to A/\theta$.
 
 By \cite[Prop 1.4.12]{Gorbunov}, a class $\mathbf{K}$ is a quasivariety iff $Co_{\mathbf{K}}(A) \subseteq Co(A)$ is an algebraic sublattice for every structure $A$. Thus the question is, whether the inclusion $Co_{\mathbf{K}}(A) \subseteq Co(A)$ preserves directed suprema (which could be used to show closure of $\bf K$ under directed colimits). A related question is whether the compact elements of $Co_{\bf K}(A)$ are exactly the finitely generated ${\bf K}$-congruences, i.e. the images under the lattice reflection of the compact elements of $Co_{\Sigma\text{-Str}}(A)$.
\end{Obs}

\subsection{Congruence filter pairs}

\begin{Df}
 Let $\Sigma$ be a signature and ${\bf K}$ a class of $\Sigma$-structures such that the association $A \mapsto Co_{\bf K}(A)$ is (the object part of) a functor from $\Sigma$-structures to algebraic lattices. A filter pair of the form $(Co_{\bf K}(A), i)$ is called congruence filter pair.
\end{Df}

In all examples in this article ${\bf K}$ will actually be a quasivariety, and in this case the set $Co_{\bf K}(A)$ of congruences relative to ${\bf K}$, ordered by inclusion, is an algebraic lattice. The next proposition shows how in this case a presentation by a congruence filter pair is closely linked to the algebraizability of a logic.

\begin{Prop}\label{2.3}\label{CriterionAlgebraizable}
Let $\Sigma$ be a signature and ${\bf K}\subseteq\Sigma\text{-Str}$ a quasivariety. If $(Co_{{\bf K}},i)$ is a filter pair, such that $i_{\Fm{}{X}}\colon Co_{\bf K}(\Fm{}{X}) \to \wp(\Fm{}{X})$ is injective, then the associated logic is algebraizable.
\end{Prop}

\Dem
We know from Lemma \ref{Lemi-FiltersAreTheImageOfi} that $i^{{\bf K}}_{Fm}[Co_{{\bf K}}(Fm)]=Th(l_{{\bf K}})$. As $i^{{\bf K}}_{Fm}$ is injective, we have that $i^{{\bf K}}_{Fm}$ is bijective. Then $i^{{\bf K}}_{Fm}$ is an isomorphism. Now let $\sigma\in hom(Fm,Fm)$. As $i^{{\bf K}}$ is a natural transformation we have the following commutative diagram:
\[
\xymatrix{
Fm\ar[d]_{\sigma}&Co_{{\bf K}}(Fm)\ar[r]^{i^{{\bf K}}_{Fm}}&(\mathcal{P}(Fm);\subseteq)\\
Fm&Co_{{\bf K}}(Fm)\ar[u]^{Co_{{\bf K}}(\sigma)}\ar[r]_{i^{{\bf K}}_{Fm}}&(\mathcal{P}(Fm);\subseteq)\ar[u]_{\sigma^{-1}}
}
\]
Note that $\sigma^{-1}(T)\in Th(l_{{\bf K}})$ for any $T\in Th(l_{k})$. Therefore $i^{{\bf K}}_{Fm}$ is a isomorphism such that commutes with substitution. By the isomorphism theorem \ref{IT1}, $l_{{\bf K}}$ is an algebraizable logic.\qed

\begin{Obs}
The above proposition gives us an alternative proof of theorem 5.2 \cite{BR}. The condition of injectivity assumed here is exactly the condition assumed in loc. cit. to get algebraizability.
\end{Obs}

%
%

The next theorem will provide us with an ample supply of examples of congruence filter pairs.

\begin{Teo}\label{equations}\label{TheoremLogicsFromEquations}
Let $\Sigma$ be a signature, ${\bf K}\subseteq\Sigma\text{-Str}$ a quasivariety and $\tau$ a finite set of equations in at most one variable. The map $i^{{\bf K}}=(i^{{\bf K}}_{M})_{M\in\Sigma\text{-Str}}$ where:
\[\begin{array}{rcl}
i^{\tau}_{M}:Co_{{\bf K}}(M)&\to&(\mathcal{P}(M),\subseteq)\\
\theta&\mapsto&\{m\in M;\ \tau^{M}(m)\in\theta\}
\end{array}\]
is a natural transformation and for any $M\in \Sigma\text{-Str}$, $i^{\tau}_{M}$ preserves arbitrary infima and directed suprema, i.e., $(Co_{{\bf K}},i^{\tau})$ is a filter pair.
\end{Teo}

\Dem
Let $f\in hom(M,N)$.
Denote here $f(\tau^{M}(m))=\{\langle f(\varepsilon^{M}(m)),f(\delta^{M}(m))\rangle;\ \langle\varepsilon,\delta\rangle\in\tau\}$.

\[
\xymatrix{
M\ar[d]_{f}&Co_{{\bf K}}(M)\ar[r]^{i^{{\bf K}}_{M}}&(\mathcal{P}(M);\subseteq)\\
N&Co_{{\bf K}}(N)\ar[u]^{Co_{{\bf K}}(f)}\ar[r]_{i^{{\bf K}}_{N}}&(\mathcal{P}(N);\subseteq)\ar[u]_{f^{-1}}
}
\]

For $\theta\in Co_{{\bf K}}(N)$ we have

$$\begin{array}{rcl}
f^{-1}(i^{{\bf K}}_{N}(\theta))&=&f^{-1}(\{n\in N;\ \tau^{N}(n)\subseteq\theta\})\\
&=&\{m\in M;\ \tau^{N}(f(m))\subseteq\theta\}\\
&=&\{m\in M;\ f(\tau^{M}(m))\subseteq\theta\}\\
&=&\{m\in M;\ \tau^{M}(m)\subseteq Co_{{\bf K}}(f)(\theta)\}\\
&=& i^{{\bf K}}_{M}(f^{-1}(\theta))
\end{array}$$

Thus $i$ is a natural transformation.

Consider a family of congruences $\theta_i \in Co_{{\bf K}}(M)\ \ (i \in I)$.

$$\begin{array}{rcl}
i^{{\bf K}}_{M}(\bigcap_{i\in I}\theta_i)&=&\{m\in M;\ \tau^{M}(m)\subseteq\bigcap_{i\in I}\theta_i\}\\
&=&\{m\in M;\ \forall i \in I: \tau^{M}(m)\subseteq \theta_i\}\\
&=&\bigcap_{i\in I} \{m\in M;\ \tau^{M}(m)\subseteq \theta_i\}
\end{array}$$ 

Thus $i_{M}$ preserves arbitrary infima.

Now let $U=\{\theta_{i};\ i\in I\}$ be an upwards directed set.

$$\begin{array}{rcl}
i^{{\bf K}}_{M}(\bigvee U)&=&\{m\in M;\ \tau^{M}(m)\subseteq \bigvee U\}\\
&=&\{m\in M;\ \tau^{M}(m)\subseteq \bigcup_{i\in I}\theta_{i}\}\\
&=&\bigcup_{i\in I}\{m\in M;\ \tau^{M}(m)\subseteq\theta_{i}\}
\end{array}$$
Here the last equality holds because $U$ is a directed set and the set of equations $\tau$ is finite, hence if each equation in $\tau(m)$ is contained in some $\theta_i$ there is a $\theta$ containing all those $\theta_i$ and hence all equations. 
Thus $i_{M}$ preserves directed suprema. \qed


\begin{Df}\label{DefinitionEquationalFilterPair}
A filter pair $(G,i)$ arising as in Theorem \ref{TheoremLogicsFromEquations} will be called an equational filter pair. The transformation $i$ will be said to be given by the equations $\tau$. We will sometimes let the defining equations be part of the notation for an equational filter pair, as in $(Co_{\bf K}, i^\tau)$ or $(Co_{\bf K}, i^{\delta \approx \epsilon})$
\end{Df}

In Thm. \ref{TeoEveryLogicWithAlgebraicSemanticsComesFromEqunlFilterPair} below we will show that the logics admitting a presentation by a congruence filter pair, are exactly those admitting an algebraic semantics. For example, every algebraizable logic admits such a presentation, with the equations given by the algebraizig pair.

For concreteness, here is a typical example of a logic with a presentation by a congruence filter pair.

\begin{Ex}\label{ExampleImplicationlessFragmentOfIPC}
Consider the signature $\Sigma = \{\wedge, \vee, \neg, \top, \bot \}$ and let ${\bf K}$ be the variety of pseudocomplemented distributive lattices. Using the single equation $\langle x, \top \rangle$ 

By \cite[Thm. 2.6]{BP1}, the logic associated to the filter pair $(Co_{{\bf K}},i)$ is the implicationless fragment of intuitionistic propositional logic $\mathbf{IPL^*}$. 
In \cite[Thm. 5.13]{BP1} it is shown that this logic is not protoalgebraic.
\end{Ex}

%
%
%
%
%


One further, admittedly artifical, example of a logic that is neither protoalgebraic nor truth-equational but has an algebraic semantics:

\begin{Ex}\label{ExampleSuccessoralgebra}
Consider the signature $\Sigma=(\Sigma)_{n\in\omega}$ which $\Sigma_{1}=\{s\}$ and $\Sigma_{n}=\emptyset$ for all $n\neq 1$. A $\Sigma$-algebra is simply a set with an endomorphism. The free $\Sigma$-algebra $F_{\Sigma}(X)$  on countably many generators $X=\{x_0, x_1, \ldots\}$ is isomorphic to $\mathbb{N}^X$, the disjoint union of countably many copies of the natural numbers, where $x_i \in X$ corresponds to $0$ in the $i$th copy, and the endomorphism $s$ acts as the successor function on each copy.

Consider the functor $Co:=Co_{\Sigma-Alg}$ and the map $i:Co \Rightarrow (\mathcal{P}(\ ),\subseteq)$ given as above by $\tau=\{\langle x,s(x)\rangle\}$. We denote the logic associated to this filter pair by $l_s$.

For a congruence $\theta$ on a $\Sigma$-algebra $A$ one has $i(\theta)=\{ a \in A \mid \bar{a}=s(\bar{a}) \text{ in } A/\theta \}$. Thus $i$ assigns to a congruence $\theta$ the set of elements of $A$ which become $s$-fixed points in $A/\theta$.

The logic $l_s$ is not protoalgebraic because it has no theorems: the theorems of $l_s$ are exactly the elements of $i(\theta_{min})$, where $\theta_{min}$ is the minimal congruence relation on $F_{\Sigma}(X)$. As $F_{\Sigma}(X)/\theta_{min} = F_{\Sigma}(X)$, this is exactly the set of $s$-fixed points in $F_{\Sigma}(X)$, which is empty.
\end{Ex}

After having seen, in Prop. \ref{CriterionAlgebraizable}, a criterion for the algebraizability of a logic in terms of a presentation by an equational filter pair, in the next few items we discuss such a criterion for truth-equationality.

\begin{Prop}\label{truth}\label{CriterionTruthEquational}
Let $(Co_{\bf K},i)$ be an equational filter pair as in Theorem \ref{TheoremLogicsFromEquations}. If $i_A$ is surjective onto the collection of filters for every algebra $A$, then the associated logic is truth-equational.
\end{Prop}
\Dem
By Definition \ref{SpecialLogics} a logic is truth-equational if and only if there is a set of equations defining the filters of all its reduced matrices. But if $i$ is surjective, then \emph{all} filters of all matrices are given by the equations defining $i$.
\qed

The following example shows that the condition of $i$ being surjective onto filters is not always satisfied, and moreover, that logics associated to equational filter pairs need not be truth-equational.

\begin{Ex}\label{ExampleSuccessorLogicIsNotTruthEquational}
Consider the logic $l_s$ of example \ref{ExampleSuccessoralgebra}. Take the $\Sigma$-algebra $A=\mathbb{N} \cup \{z\}$, where $z$ is an $s$-fixed point and $\mathbb{N}$ is endowed with the (fixed point-free) successor operation. Then $\{z\}$ is a filter on $A$ (belonging to the minimal congruence relation, whose quotient is just $A$ itself). Now $\Omega^A(\{z\})$ is the coarsest congruence relation which does not relate $z$ with any other element. This is the congruence relation identifying all the elements of $\mathbb{N}$ and leaving $z$ alone. Therefore $A/\Omega^A(\{z\})$ has two elements (both of which are $s$-fixed points). Note that $i_A (\Omega^A(\{z\}))$ has two elements and is not $\{z\}$, so that the natural candidate $\tau$ for a set of truth-equations does not work.

There is also no other set of equations defining this filter in the way required for truth-equationality: any equation satisfied by the element $\bar{z}$ in $A/\Omega^A(\{z\})$ is also satisfied by the other element, because the permutation of the two elements is an automorphism of the algebra $A/\Omega^A(\{z\})$.
\end{Ex}

There is a natural and large class of examples where the condition of Prop. \ref{CriterionTruthEquational} is satisfied:

\begin{Prop}\label{PointedQuasiVarOmegaIsRetraction}
Let ${\bf K}$ be a pointed quasivariety, i.e. a quasivariety over a signature with a constant symbol $0$. Consider the set of equations $\tau=\{\langle x,0\rangle\}$. Then $\Omega^{A}$ is a right inverse of $i^{\tau}_{A}$ for any $A\in\Sigma\text{-Str}$. In particular $i_A$ is surjective onto the collection of filters for every $\Sigma$-algebra $A$.
\end{Prop}
\Dem
Note that in the associated logic $l$ we have that $\vdash_l0$, $x,\varphi(x,\bar{z})\vdash_l\varphi(0,\bar{z})$ and $x,\varphi(0,\bar{z})\vdash_l\varphi(x,\bar{z})$ for any $\varphi(x,\bar{z})\in Fm$. Indeed, let $\theta\in Co_{{\bf K}}(Fm)$, then $\langle 0,0\rangle\in\theta$, thus $0\in i^{\tau}_{Fm}(\theta)$ and thus $\vdash_l 0$. Now let $\theta\in Co_{{\bf K}}(Fm)$ and suppose that $x,\varphi(x,\bar{z})\in i^{\tau}_{Fm}(\theta)$, then $\langle x,0\rangle\in\theta$ and $\langle \varphi(x,\bar{z}),0\rangle\in\theta$. Since $\theta$ is a congruence, we have that $\langle\varphi(x,\bar{z}),\varphi(0,\bar{z})\rangle\in\theta$. Therefore $\langle\varphi(0,\bar{z}),0\rangle\in\theta$, so $\varphi(0,\bar{z})\in i^{\tau}_{Fm}(\theta)$. Hence $x,\varphi(x,\bar{z})\vdash_l\varphi(0,\bar{z})$. The same proof can be used to prove that $x,\varphi(0,\bar{z})\vdash_l\varphi(x,\bar{z})$.

Now we are able to prove that for any $A\in\Sigma\text{-Str}$ and $F\in Fi_{l_{{\bf K}}}(A)$, $F=i^{\tau}_{A}(\Omega^{A}(F))$. Let $a\in F$ and $\varphi(x,\bar{z})\in Fm$. Let $\bar{c}\in A$ and suppose that $\varphi^{A}(a,\bar{c})\in F$. Since $x,\varphi(x,\bar{z})\vdash_l\varphi(0,\bar{z})$, we have that $\varphi^{A}(0,\bar{c})\in F$. Analogously we have that if $\varphi^{A}(0,\bar{c})\in F$, $\varphi^{A}(a,\bar{c})\in F$. Hence $\langle a,0\rangle\in\Omega^{A}(F)$. By definition of $i^{\tau}_{A}$ we have $a\in i^{\tau}_{A}(\Omega^{A}(F))$. Thus $F\subseteq i^{\tau}_{A}(\Omega^{A}(F))$. Let $a\in i^{\tau}_{A}(\Omega^{A}(F))$, then $\langle a,0\rangle\in\Omega^{A}(F)$. Since $\vdash_l 0$, we have $0\in F$, therefore $a\in F$. Hence $i^{\tau}_{A}(\Omega^{A}(F))\subseteq F$. \qed

With this we have arrived at the well-known fact that assertional logics are truth-equational:

\begin{Cor}\label{CriterionTruthEquationalPointedQuasiVar}
Let ${\bf K}\subseteq\Sigma\text{-Str}$ be a pointed quasivariety. Then the logic associated to the equation $\tau=\{\langle x,0\rangle\}$ (as in Theorem \ref{TheoremLogicsFromEquations}) is truth-equational.
\end{Cor}
\Dem Follows from Prop. \ref{PointedQuasiVarOmegaIsRetraction} and Prop. \ref{truth}.\qed

\section{Filters and adjunctions for equational filter pairs}\label{filteradjunctionequationalfilterpairs}


This section provides  a closer analysis of the left adjoint $\Xi_A$ of $i_A$ for equational filter pairs, its relation to the Leibniz operator and some consequences for the filters of the associated logics.

The main statements are Theorems \ref{TheoremInclusionsXiIdLeibniz} and \ref{ThmRegWeakAlgebraizableIffOmegaLeftAdjointOfI}.







\begin{Prop}\label{PropLogicsFromEquationalFilterPairsAreFilterWeakEquivalential}
 For a logic presented by an equational filter pair $(Co_{\bf K}, i^{\tau})$, a $\Sigma$-structure $A$ and a congruence $\theta \in Co_{\bf K}(A)$, the filter $i^{\tau}_A(\theta)$ is compatible with the congruence relation $\theta$ (in the sense of Def. \ref{DefCompatibleCongruence}).
\end{Prop}
\begin{proof}
 Let the filter pair be given by the equations $\tau = \langle \epsilon(x), \delta (x) \rangle$. Let $\theta \in Co_{\bf K}(A)$ and $\langle \varphi,\psi \rangle \in \theta$. Then we have 
$$\begin{array}{rcl}
 \varphi \in i(\theta)& \Leftrightarrow &  \langle \epsilon(\varphi), \delta(\varphi) \rangle \in \theta \text{ (by definition of $i(\theta)$)}\\
& \Leftrightarrow &  \epsilon(\bar{\varphi})=\delta(\bar{\varphi}) \in A/\theta \\
& \Leftrightarrow &  \epsilon(\bar{\psi})=\delta(\bar{\psi}) \in A/\theta  \text{ (because $\bar{\varphi} = \bar{\psi} \in A/\theta$)} \\
& \Leftrightarrow & \psi \in i(\theta) \text{ (by definition of $i(\theta)$)} \\
\end{array}$$
 \end{proof}

 Next we give an explicit description of the operator  $\Xi_A$ of Def. \ref{DefLeftAdjointOfAFilterPair} in the case of an equational filter pair.

\begin{Prop}\label{PropFirstFormulaForLeftAdjointForEquationalFilterPairs}
Let $(Co_{\bf K},i)$ be an equational filter pair, with $i$ given by the equations $\delta(x)=\epsilon(x)$. Then for a $\Sigma$-algebra $A$ the left adjoint $\Xi_A$ to $i_A$ maps $S \subseteq A$ to the congruence relation generated by $\{\langle \delta(s), \epsilon(s) \rangle \mid s \in S\}$.
\end{Prop}
\Dem
This follows immediately from the general description of Theorem \ref{TheoremExistenceOfLeftAdjoint} of left adjoints between complete partially ordered sets:

\[
\begin{array}{rcl}
\Xi_A(S) &=& \bigwedge \{ \theta \in Co_{\bf K}(A) \mid S \subseteq i_A(\theta) \} \\ 
&=& \bigwedge \{ \theta \in Co_{\bf K}(A) \mid S \subseteq \{a \in A \mid \delta(\bar{a})=\epsilon(\bar{a}) \text{ in } A/\theta \} \} \\
&=& \bigwedge \{ \theta \in Co_{\bf K}(A) \mid \{ \langle \delta(s), \epsilon(s)\rangle \mid s \in S \}   \subseteq  \theta\} 
\end{array}
\]
\qed

Recall from Def. \ref{DefEquationalConsequence} the notation $Cn_{\bf K}(R)$ for the set of equations that is implied in the quasivariety ${\bf K}$ by a set of equations $R$. Also recall from and Lemma \ref{LemmaEquationalClosureEqualsGeneratedCongruence} that, if one identifies equations with pairs of elements, $Cn_{\bf K}(R)$ is exactly the  congruence relative ${\bf K}$ generated by the set of pairs $R$. 
Blok-Pigozzi, in \cite[p.28]{BP1}, on the formula algebra $\Fm{}{X}$ define the operator $$\Omega_{\bf K} \colon \wp(\Fm{}{X}) \to Co_{\bf K}(\Fm{}{X}), \ \ \ \ T \mapsto Cn_{\bf K}(\{ \delta(a)\approx\epsilon(a) \mid a \in T \}).$$ This makes sense on any algebra and gives us an alternative description of $\Xi_A$:

\begin{Cor}\label{PropOurLeftAdjointIsBlokPigozzisOmegaK}
For an algebra $A$ and a subset $T \subseteq A$ one has $\Xi_A (T) = Cn_{\bf K}(\{ \delta(a)\approx\epsilon(a) \mid a \in T \})$.
\end{Cor}
\Dem
Follows immediately from Lemma \ref{LemmaEquationalClosureEqualsGeneratedCongruence} and Proposition \ref{PropFirstFormulaForLeftAdjointForEquationalFilterPairs}.
\qed

With Corollary \ref{PropOurLeftAdjointIsBlokPigozzisOmegaK} we can now pin down the class of logics admitting a presentation by an equational filter pair.

\begin{Teo}\label{TeoEveryLogicWithAlgebraicSemanticsComesFromEqunlFilterPair}
Let $\Sigma$ be a signature, $l$ a logic over $\Sigma$, ${\bf K}$ a class of $\Sigma$-algebras and ${\bf K}^Q$ the quasivariety generated by ${\bf K}$.
 Then $l$ has an algebraic semantics in ${\bf K}$, if and only if it has a presentation by an equational filter pair $(Co_{{\bf K}^Q},i)$.
\end{Teo}
\Dem
Assume that $l$ has an algebraic semantics in ${\bf K}$ given by a set $\{\delta(x)=\epsilon(x)\}$ of equations over $\Sigma$. 

Using Theorem \ref{TheoremLogicsFromEquations} we can define the filter pair $(Co_{{\bf K}^Q}, i^{\delta = \epsilon})$ and compare its associated logic with $l$. Let $\Xi_{\Fm{}{X}}$ be the left adjoint of $i_{\Fm{}{X}}$.

By \cite[Cor 2.3]{BP1} one can replace ${\bf K}$ by ${\bf K}^Q$, hence

 $$\begin{array}{rclr}
\Gamma \vdash_l \varphi  & \Leftrightarrow  & \{\delta(\gamma) = \epsilon(\gamma) \mid \gamma \in \Gamma\} \vDash_{{\bf K}} \delta(\varphi) = \epsilon(\varphi) &   \\
 & \Leftrightarrow  & \{\delta(\gamma) = \epsilon(\gamma) \mid \gamma \in \Gamma\} \vDash_{{\bf K}^Q} \delta(\varphi) = \epsilon(\varphi)&   \text{\cite[Cor 2.3]{BP1}} \\
 & \Leftrightarrow  & \delta(\varphi) = \epsilon(\varphi) \in Cn_{{\bf K}^Q}(\{\delta(\gamma) = \epsilon(\gamma) \mid \gamma \in \Gamma\})&   \text{def. of $Cn_{{\bf K}^Q}$} \\
 & \Leftrightarrow  & \delta(\varphi) = \epsilon(\varphi) \in \Xi_{\Fm{}{X}}(\Gamma)&   \text{Cor. \ref{PropOurLeftAdjointIsBlokPigozzisOmegaK}} \\
 & \Leftrightarrow  & \varphi \in i^{\delta = \epsilon}(\Xi_{\Fm{}{X}}(\Gamma))&    \text{def. of $i^{\delta = \epsilon}$} \\
 & \Leftrightarrow  & \Gamma \vdash_{(Co_{{\bf K}^Q}, i^{\delta = \epsilon})} \varphi &   \text{Prop. \ref{PropDescriptionOfConsequenceRelationByClosureOperator}}
\end{array}$$

Vice versa it follows from the last four equivalences that a logic presented by an equational filter pair has an algebraic semantics.
\qed

We dedicate the rest of the section to a study of the properties of the family of adjunctions $i^\tau, \Xi$. To begin, many of the results of Chapter 3 of \cite{BP1} on the operator $\Omega_{\bf K}$, there proven directly from the definition and only for the formula algebra $\Fm{}{X}$, follow now immediately from the fact that it is left adjoint to $i^\tau$.

\begin{Lem}\label{LemmaConsequencesOfAdjointnessForEquationalFilterPairs}
Let $(Co_{\bf K}, i^{\tau})$ be an equational filter pair and $A$ a $\Sigma$-structure. 
\begin{enumerate}[(i)]
  \item \emph{\cite[Lemma 3.3(i)]{BP1}} Write $C := i_A^\tau \circ \Xi_A$ for the closure operator of the abstract logic on $A$. Then for $\Gamma \subseteq A$ we have $\Xi_A(C(\Gamma)) = \Xi_A(\Gamma)$\footnote{Blok-Pigozzi's original statement follows together with Prop. \ref{PropDescriptionOfConsequenceRelationByClosureOperator}}.
  
  \item \emph{\cite[Lemma 3.3(ii)]{BP1}} The map $\Xi_A$ preserves arbitrary suprema.
  
  \item \emph{\cite[Lemma 3.3(iii)]{BP1}} The map $\Xi_A$ preserves unions of directed sets of theories of the abstract logic on $A$.
  
  \item \emph{\cite[Lemma 3.4(i)]{BP1}} For a theory $T$ of the abstract logic on $A$ one has $(i^\tau_A \circ \Xi_A )(T) = T$
  
  \item \emph{\cite[Lemma 3.4(ii)]{BP1}} For every $\theta \in Co_{\bf K}(A)$ one has $(\Xi_A \circ i^\tau)(\theta) \subseteq \theta$. Equality holds if and only if $\theta = \Xi_A(T)$ for some $A$-theory $T$.
  
  \item \emph{\cite[Lemma 3.5(i)]{BP1}} $\Xi_A$ maps the lattice of $A$-theories isomorphically to a compact and join-complete sublattice of $Co_{\bf K}(A)$.
 \end{enumerate}
\end{Lem}
\Dem
\begin{enumerate}[(i)]
 \item This is just the triangle equality $\Xi_A \circ i^\tau \circ \Xi_A = \Xi_A$ coming from the adjointness.
 
 \item The map $\Xi_A$ is a left adjoint, hence preserves arbitrary colimits.
 
 \item Since both $\Xi_A$ and $i_A^\tau$ preserve directed suprema, the composite $i_A^\tau \circ \Xi_A$ preserves directed suprema in $\wp(A)$, i.e. directed unions, and takes values in the sublattice of theories. Thus directed suprema in the lattice of theories coincide with set-theoretic unions, hence this follows from part (ii).

 \item By Lemma \ref{Lemi-FiltersAreTheImageOfi} theories of the abstract logic on $A$ (there called $i$-filters, see Def. \ref{DefiFilter}) are exactly the sets in the image of $i^\tau_A$. Hence the claim follows from the adjointness equality $i^\tau_A \circ \Xi_A \circ i^\tau = i^\tau$.
 
 \item The inclusion is the counit of the adjunction, the condition for equality is clear from the adjunction properties.

 \item Adjunctions induce isomorphisms between the sublattices occurring as the images. The map $\Xi_A \circ i^\tau$ is a coreflection onto the respective sublattice of $Co_{\bf K}(A)$ and coreflective sublattices are closed under suprema.
 \end{enumerate}
\qed

As mentioned earlier, the collection of adjoint maps $(\Xi_A)_{A \in \Sigma\textrm{-Str}}$ need not be a natural transformation with respect to preimage maps in general, and its behaviour with respect to homomorphisms of $\Sigma$-structures does not follow directly from adjointness. It is, however, of significance in the consideration of Section \ref{SectionCraigInterpolation}. The following lemma records the behaviour with respect to substitutions on the formula algebras.

\begin{Lem}\label{LemmaSubstitutionInvarianceForTheEquationalOperators}
Let $\sigma \colon \Fm{}{X} \to \Fm{}{Y}$  be a homomorphism of free $\Sigma$-algebras.
\begin{enumerate}[(i)]
 \item\label{substitutionInvarianceKConsequence} For $\Gamma \subseteq \Fm{}{X} \times \Fm{}{X}$ we have $(\sigma \times \sigma)(Cn_{\bf K}^{\Fm{}{X}}(\Gamma)) \subseteq Cn_{\bf K}^{\Fm{}{Y}}(\sigma(\Gamma))$.
 
 \item\label{equalityKConsequenceClosuresAndSubstitution} For $\Gamma \subseteq \Fm{}{X} \times \Fm{}{X}$ we have $Cn_{\bf K}^{\Fm{}{Y}}((\sigma \times \sigma)(Cn_{\bf K}^{\Fm{}{X}}(\Gamma))) = Cn_{\bf K}^{\Fm{}{Y}}(\sigma(\Gamma))$.
 
 \item\label{substitutionInvarianceDirectImageXi}  For any homomorphism of free $\Sigma$-algebras $\sigma \colon \Fm{}{X} \to \Fm{}{Y}$ and $T \subseteq \Fm{}{X}$ we have \linebreak $\Xi_{\Fm{}{Y}}(\sigma(T)) = Cn_{\bf K}^{\Fm{}{Y}}((\sigma\times \sigma)(\Xi_{\Fm{}{X}}(T)))$
 
 \item\label{substitutionInvarianceInverseImageXi}  If $\sigma \colon \Fm{}{X} \to \Fm{}{Y}$ is the homomorphism induced by an injective map of variables $X \hookrightarrow Y$ and $T \subseteq \Fm{}{Y}$, we have $\Xi_{\Fm{}{X}}(\sigma^{-1}(T)) \subseteq (\sigma\times \sigma)^{-1}(\Xi_{\Fm{}{Y}}(T))$

\end{enumerate}
 \end{Lem}

\Dem
\begin{enumerate}[(i)]
 \item Equational consequence satisfies, for any substitution $\sigma \colon \Fm{}{X} \to \Fm{}{Y}$, and set $\Gamma \cup \{\tau \} \subseteq \Fm{}{X}$ of equations, $\Gamma \vDash_{\bf K} \tau \Rightarrow \sigma\Gamma \vDash_{\bf K} \sigma \tau$. Hence we have $(\sigma \times \sigma)(Cn_{\bf K}^{\Fm{}{X}}(\Gamma)) = \{(\sigma \times \sigma)(\tau)  \mid \Gamma \vDash_{\bf K} \tau  \} \subseteq \{  \tau \mid (\sigma \times \sigma)(\Gamma) \vDash_{\bf K} \tau  \} = Cn_{\bf K}^{\Fm{}{Y}}(\sigma(\Gamma))$.
 
  \item Applying $Cn_{\bf K}^{\Fm{}{Y}}$ to the inclusion from (\ref{substitutionInvarianceKConsequence}) yields $$Cn_{\bf K}^{\Fm{}{Y}}((\sigma \times \sigma)(Cn_{\bf K}^{\Fm{}{X}}(\Gamma)))  \subseteq  Cn_{\bf K}^{\Fm{}{Y}} Cn_{\bf K}^{\Fm{}{Y}}(\sigma \times \sigma(\Gamma)) = Cn_{\bf K}^{\Fm{}{Y}}(\sigma \times \sigma(\Gamma)).$$
  On the other hand $Cn_{\bf K}^{\Fm{}{X}}(\Gamma) \supseteq \Gamma$ implies $(\sigma \times \sigma)(Cn_{\bf K}^{\Fm{}{X}}(\Gamma)) \supseteq (\sigma \times \sigma)(\Gamma)$, which implies $Cn_{\bf K}^{\Fm{}{Y}}((\sigma \times \sigma)(Cn_{\bf K}^{\Fm{}{X}}(\Gamma))) \supseteq  Cn_{\bf K}^{\Fm{}{Y}}((\sigma \times \sigma)(\Gamma)).$
 
 \item See \cite[Lemma 3.6]{BP1} for the case $X=Y$. The general proof is no different.
 
 \item  In the special case we are considering now, $\sigma^{-1}$ (resp. $(\sigma \times \sigma)^{-1}$) is just restriction to formulas (resp. pairs of formulas) in $\Fm{}{X}$, i.e. $\sigma^{-1}(T) = T \cap \Fm{}{X}$.
 
 Recall, from the text after Lemma \ref{LemmaInclusionOfVariablesGivesConservativeTranslationOfEquationalConsequences}, the notation $Cn_{\bf K}^X$ for the closure operator associated to the equational consequence relation $\vDash_{\bf K}^X$ of Def. \ref{DefEquationalConsequence}. Expressed in terms of these closure operators, Lemma \ref{LemmaInclusionOfVariablesGivesConservativeTranslationOfEquationalConsequences} says that $Cn_{\bf K}^Y(T \cap \Fm{}{X}^2) \cap \Fm{}{X}^2 = Cn_{\bf K}^X(T \cap \Fm{}{X}^2)$. Hence 
 $$Cn_{\bf K}^Y(T) \cap \Fm{}{X}^2 \supseteq Cn_{\bf K}^Y(T \cap \Fm{}{X}^2) \cap \Fm{}{X}^2 = Cn_{\bf K}^X(T \cap \Fm{}{X}^2)$$
which implies the claim.
\end{enumerate}
\qed

While in the usual treatments of abstract algebraic logic the Leibniz operator plays a prominent role, we haven't made any use of it so far. Even algebraizable logics can be characterized purely in terms of the natural transformations $i$ and $\Xi$ (namely if $\Xi(i(\theta))$ always holds, Thm. \ref{TheoremInclusionsXiIdLeibniz} a)$\Leftrightarrow$e) below).

We will now start considering the interplay of our operators with the Leibniz operator $\Omega$.

\begin{Teo}\label{TheoremInclusionsXiIdLeibniz}
Let $(Co_{\bf K},i)$ be an equational filter pair, with $i$ given by equations $\delta(x)= \epsilon(x)$. Let $A$ be a $\Sigma$-algebra and $\theta \in Co_{\bf K}(A)$ a congruence on $A$. Then $\Xi_A(i(\theta)) \subseteq \theta \subseteq \Omega_A(i_A(\theta))$. Moreover the following are equivalent:
\begin{enumerate}[(a)]
 \item one of the inclusions is an equality for every $\theta \in Co_{\bf K}(\Fm{}{X})$ (where $X$ is any enumerable set)

 \item both inclusions are equalities for all $\Sigma$-structures $A$ and every $\theta \in Co_{\bf K}(A)$
 
 \item $\Xi_A = \Omega_A$ for all $\Sigma$-structures $A$, as maps from the lattice of filters to $Co_{\bf K}(A)$
 
 \item $\Xi_{\Fm{}{X}} = \Omega_{\Fm{}{X}}$, as maps from the lattice of theories to $Co_{\bf K}(\Fm{}{X})$

 \item the logic presented by $(Co_{\bf K},i)$ is algebraizable with $i_{\Fm{}{X}}$ being a substitution preserving lattice isomorphism. 
 \end{enumerate}
\end{Teo}
\Dem
The first inclusion has been noted in Lemma \ref{LemmaConsequencesOfAdjointnessForEquationalFilterPairs}(v).

For the second inclusion recall that by Prop. \ref{PropLogicsFromEquationalFilterPairsAreFilterWeakEquivalential} $\theta$ is compatible with $i(\theta)$. But $\Omega_A(i_A(\theta))$ is the coarsest congruence relation that is compatible with $i(\theta)$, so $\theta \subseteq \Omega_A(i_A(\theta))$.

For the equivalences we prove the cycle $(a) \Rightarrow (e) \Rightarrow (c) \Rightarrow (b) \Rightarrow (d) \Rightarrow (a)$.

$(a) \Rightarrow (e):$ If one of the inclusions $\theta \subseteq \Omega(i(\theta))$, $\theta \subseteq \Xi(i(\theta))$ is an equality, then $i$ has a left inverse and hence is injective. Since by definition of the associated logic, $i$ is also surjective onto the theories, we have a substitution preserving isomorphism between the lattice of congruences and the lattice of theories. By Thm. \ref{TheoremSecondIsoThmInclUniquenessOfTheQuasivariety} the associated logic is algebraizable.

$(e) \Rightarrow (c):$ By Thm. \ref{TheoremSecondIsoThmInclUniquenessOfTheQuasivariety} for any $\Sigma$-structure $A$ the equations defining $i_A$ induce isomorphisms, and the inverse is given by the Leibniz operator $\Omega_A$. Being an inverse of $i_A$, $\Omega_A$ is in particular a left adjoint. By uniqueness of adjoints, it follows that $\Xi_A = \Omega_A$.

$(c) \Rightarrow (b):$ Clear, using the two inclusions and the equality $\Xi_A = \Omega_A$.

$(b) \Rightarrow (d):$ By hypothesis $(b)$ the maps $\Xi_{\Fm{}{X}}$ and $\Omega_{\Fm{}{X}}$ coincide on the filters in the image of $i_{\Fm{}{X}}$. But by Prop. \ref{FreeAlgebraResult} all theories are in the image of $i_{\Fm{}{X}}$, hence $\Xi_{\Fm{}{X}} = \Omega_{\Fm{}{X}}$.

$(d) \Rightarrow (a):$ Clear, using the two inclusions and the equality $\Xi_{\Fm{}{X}} = \Omega_{\Fm{}{X}}$.
\qed

%
%

\begin{Lem}\label{LemmaFiltersInTheImageOfiForReducedMatrices}
Let $(Co_{\bf K},i)$ be an equational filter pair, with $i$ given by equations $\delta(x)= \epsilon(x)$. Let $\langle A, F \rangle$ be a reduced matrix for the associated logic with $F=i_A(\theta)$ in the image of $i_A$, for a congruence $\theta \in Co_{\bf K}(A)$. Then $F=\{a \in A \mid \delta(a)=\epsilon(a) \}$.
\end{Lem}
\Dem
The assumption that $\langle A, F \rangle$ is a reduced matrix means $\Omega_A(F)=\Delta$ (where $\Delta$ denotes the diagonal, i.e. the minimal congruence relation). We therefore have $\Omega_A(F)=\Delta \subseteq \Xi_A(F)$, and thus by the inclusions of Theorem \ref{TheoremInclusionsXiIdLeibniz} $\Delta=\Xi_A(F)=\theta$. We obtain
\[
\begin{array}{rcl}
 F = i_A(\theta) &=& \{ a \in A \mid \delta(\bar{a})=\epsilon(\bar{a})\text{ in } A/\theta \} \\
&=& \{ a \in A \mid  \delta(a)=\epsilon(a) \text{ in } A \}
\end{array}
\] 
\qed

Our considerations about adjoints and filter pairs let us arrive at the following well-known fact about assertional logics.

\begin{Cor}\label{CorollaryDescriptionOfReducedMatricesForTruthEquationalFilterPairs}
Let $(Co_{\bf K},i)$ be an equational filter pair, with $i$ given by equations $x=\top$. Let $\langle A, F \rangle$ be a reduced matrix for the associated logic. Then $F=\{ \top \}$. 
\end{Cor}
\Dem
For a filter pair given by the equation $x=\top$ all filters are in the image of $i_A$ by Lemma \ref{PointedQuasiVarOmegaIsRetraction}. By Lemma \ref{LemmaFiltersInTheImageOfiForReducedMatrices} we have $F=\{a \in A \mid a=\top\}=\{ \top \}$. 
\qed

We know from Theorem \ref{TheoremInclusionsXiIdLeibniz} that if one of the inclusions $\Xi(i(\theta)) \subseteq \theta \subseteq \Omega(i(\theta))$ is actually an equality (and hence $i$ injective), then the filter pair in question presents an algebraizable logic. 
Applying $i$ to the first inclusion, by general properties of adjunctions, we obtain an equality $i(\Xi(i(\theta))) = i(\theta)$. In contrast we do not necessarily have the corresponding equality $i(\theta) = i(\Omega(i(\theta)))$ for the Leibniz operator.

\begin{Ex}
 Consider again the logic $l_s$ of Examples \ref{ExampleSuccessoralgebra} and \ref{ExampleSuccessorLogicIsNotTruthEquational}. In Example \ref{ExampleSuccessorLogicIsNotTruthEquational} we considered the algebra $A:= \mathbb{N} \cup \{z\}$ with the successor operation on $\mathbb{N}$ and $z$ a fixed point for the unary operation $s$. As $\{z\}$ is the $s$-fixed point set of $A$, we have $\{z\} = i(\Delta_A)$ (where $\Delta_A$ denotes the minimal congruence relation on $A$). As noted in Example \ref{ExampleSuccessorLogicIsNotTruthEquational}, $\Omega_A(\{z\})$, being the coarsest congruence not relating elements of $\{z\}$ with elements outside $\{z\}$, is the congruence collapsing all elements of $\mathbb{N}$ to a single element and leaving $z$ alone. The quotient $A/\Omega_A(\{z\})$ consists of two elements, both of which adre $s$-fixed points. Hence $i (\Omega_A(\{z\}))$, the set of elements of $A$ which become $s$-fixed points in $A/\Omega_A(\{z\})$, is all of $A$. Altogether we have $i (\Omega_A(i(\Delta_A))) = A \supsetneq \{z\} = i(\Delta_A)$.
 \end{Ex}

Continuing to consider the previous example we may apply $\Omega$ to the inequality $i (\Omega_{A}(i(\Delta_A))) = A \supsetneq \{z\} = i(\Delta_A)$ and obtain $\Omega(i (\Omega(i(\Delta_A)))) = \Omega(A) = A \times A \supsetneq \Omega(\{z\}) = \Omega(i(\Delta_A))$ where the strict inequality comes from the fact that the quotient of the left hand congruencce yields a one element set and the quotient of the right hand congruence yields a two element set. This shows that one also in general does not have equality in the inclusion $\Omega \circ i \circ \Omega \circ i \supseteq \Omega \circ i$. However, in all examples we considered the iterated application of $\Omega \circ i$ stabilizes at some point (in the above example this happens in the next step).
 
\begin{Que}
Does iterated application of $\Omega \circ i$ always stabilize, i.e. is there for any congruence $\theta$ an $n \in \mathbb{N}$ such that $(\Omega \circ i)^{\circ n}(\theta) = (\Omega \circ i)^{\circ (n+1)}(\theta)$? If yes, is there a uniform bound on $n$ for all subsets $F$ of a fixed algebra $A$?
\end{Que}
 
While the condition $\theta = \Omega(i(\theta))$ is equivalent to algebraizability, the weaker condition $i(\theta) = i(\Omega(i(\theta)))$ can, however, also hold for non-algebraizable logics: We have seen that for assertional logics we have $i(\Omega(F))=F$, hence $i = i\circ \Omega\circ i$.

 Let $(Co_{\bf K},i)$ be an equational filter pair, with $i$ given by equations $x=\top$. Then we have $i = i\circ \Omega\circ i$. 

\begin{Que}
 Can one give a meaningful characterization of the class of equational filter pairs for which $i = i\circ \Omega\circ i$ holds? Is it precisely the class of filter pairs given by equations of the form $x=\top$?
\end{Que}

The following discussion will shed some more light on this question.

We have seen that the operator $\Xi \circ i$ associates to a congruence $\theta$ the \emph{smallest} congruence in whose quotient the same elements get mapped to solutions of the equations $\tau$ set as in the quotient by $\theta$ -- in particular this smallest congruence exists. If $i$ has a right adjoint $R$, then the operator $R \circ i$ associates to a congruence $\theta$ the \emph{biggest} congruence which has the same $\tau$-solution set as $\theta$. This biggest congruence exists if and only if there is a right adjoint, but this need not be the case in general.

One may suggest the mental image that the Leibniz operator is trying to be a right adjoint of $i$, but may fail to be so. Indeed, the failure may lie in the fact that it is not even an order preserving map for non-protoalgebraic logics. For protoalgebraic logics, the Leibniz operator is not only order preserving, but also preserves arbitrary infima, and hence by Theorem \ref{TheoremExistenceOfLeftAdjoint} has a left adjoint. The description of the left adjoint offered by that same theorem resembles the description of $i$ in the case of an assertional logic, as we will see below.

\begin{Prop}
 Let $l$ be a protoalgebraic logic. Then the Leibniz operator from the \emph{filters} to the congruences has a left adjoint given by $L(\theta) = \bigcap \{ F \ \mid \ \theta \textrm{ does not relate elements of $F$ with elements outside $F$}   \}$.
\end{Prop}

\Dem
Since it preserves arbitrary infima, $\Omega$ has a left adjoint $L$. By the adjunction formula Thm. \ref{TheoremExistenceOfLeftAdjoint} it is given by:

$$\begin{array}{rcl}
    L(\theta) &=& \bigwedge \{ F \, \mid \, \theta \subseteq \Omega(F) \} \\
    &=& \bigwedge \{ F \, \mid \, \theta \subseteq \vee \{\rho \mid F \textrm{ is a union of } \rho-\textrm{equivalence classes}\} \} \\
    &=& \bigcap \{ F \ \mid \ \theta \textrm{ does not relate elements of $F$ with elements outside $F$}   \}
\end{array}$$

Arbitrary intersections of sets that are compatible with a given equivalence relation $\theta$ are compatible with $\theta$ again. Arbitrary intersections of filters are filters again, so replacing the infimum on the last step with an intersection is ok. \qed





We would like to relate this to the right adjoint $i$ occurring in a filter pair. Consider the case where $i$ is given by the equation $\langle x, \top \rangle$, then we have the formula $i(\theta) = \{a \ \mid \ a \theta \top\} = [\top]_{\theta}$. We have $i(\theta) \subseteq F$ for each filter $F$ such that $\theta \subseteq \Omega(F)$, thus $i(\theta) \subseteq L(\theta)$. The opposite inclusion, $L(\theta) \subseteq i(\theta)$, holds too, since $[\top]_{\theta}$ occurs among the filters over which we take the intersection in the formula for $L(\theta)$ (this inclusion is also a consequence of the adjunction, since $\theta \subseteq \Omega(i(\theta))$, because $\theta$ is compatible with $i(\theta)$). {\bf Thus $i=L$ !}


The hypothesis that our logic admits a presentation by a congruence filter pair with $i$ given by the equation $\langle x,\top\rangle$ is equivalent to asking that it is assertional.
By \cite[6.125]{Fon}, a logic is assertional and protoalgebraic if and only if it is regularly weakly algebraizable.

Thus we arrive at the following result:

\begin{Teo}\label{ThmRegWeakAlgebraizableIffOmegaLeftAdjointOfI}
 Let $l$ be the logic associated to a congruence filter pair $(Co_{\mathbf{K}},i)$, with $i$ given by the equation $\langle x,\top\rangle$.  Then $l$ is regularly weakly algebraizable if and only $\Omega$ is a left adjoint of $i$.
\end{Teo}

\Dem
One direction has been shown in the above discussion: If a logic is regularly weakly algebraizable, then it is protoalgebraic and we have seen above that then the left adjoint of $\Omega$ coincides with $i$.

In the other direction note that if $i$ is left adjoint to $\Omega$, then $\Omega$, being a right adjoint, preserves arbitrary infima and hence the logic is protoalgebraic. Assertionality follows directly from the filter pair being given by the equation $\langle x,\top\rangle$.
\qed

Note that the equation $i = i\circ \Omega\circ i$, shown independently before, now follows immediately from the adjointness of $i$ and $\Omega$.


\begin{Lem}\label{protoalgebraic}
$\Omega \colon Fi \to Co_K$ has a left adjoint if and only if the logic is protoalgebraic.
\end{Lem}
\Dem
$\Omega$ has a left adjoint if and only if it preserves arbitrary infima (by Thm. \ref{TheoremExistenceOfLeftAdjoint}) if and only if the logic is protoalgebraic, by Thm. 6.4 in \cite{Fon}.  
\qed

\begin{Prop}\label{PropProtoAlgAndReducedAlgsEqualKThenAlgebraizable}
Let $l$ be the logic associated to a filter pair $(Co_K, i)$, with $i$ given by the equation $x=\top$. Suppose that $l$ is protoalgebraic and that the class of reduced algebras coincides with $K$. Then $l$ is algebraizable.
\end{Prop}
\Dem
We know already that $l$ is weakly regularly algebraizable (by a combination of the previous Lemma and Proposition). I remains to show that $l$ is also equivalential.

For a weakly algebraizable logic, by \cite[Thm. 6.117]{Fon} the Leibniz operator is a bijection $\Omega^A\colon Fi_l A \simeq Con_{Alg^{*}\,l}A$ from filters to congruences
relative the class of reduced algebras for $l$. If this class coincides with $\mathbf{K}$, then we have that $\Omega$ is a bijection $\Omega^A\colon Fi_l A \cong Con_{\mathbf{K}}A$. Its adjoint $i$ is then also a bijection and hence the logic is algebraizable. Alternatively, since $i$ is the adjoint, and hence inverse, of $\Omega$, and $i$ is a natural transformation, $\Omega$ must be natural as well, which means that $l$ is equivalential by \cite[Cor. 6.69]{Fon}.\qed

Note that, in this case, {\bf K} is a quasivariety that is the equivalent algebraic semantics for $l$.


\begin{Ex}
 Consider the filter pair $(Co_K, i^\tau)$ where {\bf K} is the variety of all groups, axiomatized with a constant symbol $e$ for the neutral element, and $i^\tau$ given by the equation $\tau=\langle x, e \rangle$.
 
 The associated logic is protoalgebraic with implication connective $y^{-1}x$: This follows from the fact that if $y^{-1}x=e$ and $x^{-1}z=e$ then also $y^{-1}z = y^{-1}xx^{-1}z=ee = e$.  Furthermore, every group has $\{e\}$ as a reduced filter, so every group is the underlying algebra of a reduced matrix. By Prop. \ref{PropProtoAlgAndReducedAlgsEqualKThenAlgebraizable} the associated logic is algebraizable.
 
 A similar reasoning applies to the filter pair $(Co_K, i^\tau)$ with {\bf K} the variety of all rings, and $i^\tau$ given by the equation $\tau=\langle x, 0 \rangle$.
\end{Ex}



On the other hand, classical logic $l_c$ has as algebraic semantics the quasivariety generated  by the orthomodular lattice $O_6$ (see \cite[Ex. 4.79 ]{Fon}), which is not the class of boolean algebras. Thus the criterion of Prop. \ref{PropProtoAlgAndReducedAlgsEqualKThenAlgebraizable} is sufficient but not necessary for a logic to be algebraizable. This comes from the possibility of giving an algebraic semantics to an algebraizable logic that is not an equivalent algebraic semantics.

\section{Craig entailment interpolation property and filter pairs}\label{SectionCraigInterpolation}

In this section we present a correspondence between the Craig interpolation property for a logic associated to an equational filter pair, and the amalgamation property in the class ${\bf K}$ for whose congruences the filter pair is defined.

We recall the two central notions for this chapter.

\begin{Df}
\begin{itemize}\label{DefinitionMatrixEmbeddingCraigEntailmentAmalgamationProperty}
\item A logic $l$ has the \emph{Craig entailment interpolation property} if for every set of formulas $\Gamma$, with variables $var(\Gamma)$ and every formula $\varphi$ with variables $var(\varphi)$, if $\Gamma\vdash \varphi$ then there is a set of formulas $\Gamma'$ with the variables in $var(\Gamma)\cap var(\varphi)$ such that $\Gamma\vdash \Gamma'$ and $\Gamma'\vdash \varphi$.
\item We shall say that a class of algebras ${\bf K}$ has the \emph{amalgamation property} if given $A,B,C\in {\bf K}$ and injective homomorphisms $i_{B} \colon A \to B$, $i_{C} \colon A \to C$, there exist an algebra $D \in {\bf K}$ and injective homomorphisms $e_B:B\to D$, $e_C:C\to D$ such that $e_{B}\circ i_{B}=e_{C}\circ i_{C}$.
\end{itemize}
\end{Df}

We will employ the following auxiliary property which implies Craig interpolation.

\begin{Df}
  [Def. 3.4 of \cite{CzP}] A logic $l$ has the \emph{flat theory amalgamation property} if for every two non-disjoint sets of variables $X$ and $Y$, and every $l$-filter $T$ of the formula algebra $\Fm{}{X}$ there is an $l$-filter $R$ of the formula algebra $\Fm{}{X\cup Y}$ such that $R\cap \Fm{}{X}=T$ and $R\cap \Fm{}{Y}=Fi_{l}[T\cap \Fm{}{X\cap Y}]=\bigcap\{T'\in Th(\Fm{}{Y});\ T\cap \Fm{}{X\cap Y}\subseteq T'\}$
\end{Df}

\begin{Lem}[\cite{CzP}, Thm. 3.5]\label{2.17}
If a logic $l$ has the flat theory amalgamation property, it has the Craig entailment interpolation property.
\end{Lem}
\Dem Let $\Gamma\cup\{\varphi\}$ set of formulas such that $\Gamma\vdash \varphi$, $X=var(\Gamma)$ and $Y=var(\varphi)$. Consider the $l$-filter $T$ of $\Fm{}{X}$ generated by $\Gamma$. By the flat theory amalgamation property, there is an l-filter $R$ of $\Fm{}{X\cup Y}$ such that $R\cap \Fm{}{X}=T$ and $R\cap \Fm{}{Y}=Fi_{l}[T\cap \Fm{}{X\cap Y}]$. Note that $\varphi\in R\cap \Fm{}{Y}$, and consider $\Gamma'=T\cap \Fm{}{X\cap Y}$. It is clear that $\Gamma\vdash \Gamma'$ and $\Gamma'\vdash\varphi$. \qed

We now introduce a property which ensures that from a situation in which one can ask for flat theory amalgamation one can pass to a situation in which one can apply amalgamation of algebras. It is the content of Theorem \ref{TheoremMatrixAmalgamationImpliesCraigInterpolation} below that whenever this property is satisfied, the amalgamation property for a quasivariety ${\bf K}$ implies Craig interpolation for the logic presented by a filter pair $(Co_{\bf K}, i)$.

\begin{Df}\label{DefTheoryLiftingProperty}
 Let ${\bf K}$ be a class of algebras. We say that a congruence filter pair $(Co_{{\bf K}},i)$ has the \emph{theory lifting property}, if the following holds:
 
 \medskip
 
 Given sets $X,Y$ such that $Z:=X \cap Y \neq \emptyset$ and a theory $T \subseteq \Fm{}{X}$ of the logic associated to $(Co_{{\bf K}},i)$, define $T'':=T \cap \Fm{}{Z}$ and $T' \subseteq \Fm{}{Y}$ to be the smallest $\Fm{}{Y}$-theory containing $T''$. 
 
 Then there exist $\theta_T \in Co_{\bf K}(\Fm{}{X})$, $\theta_{T''} \in Co_{\bf K}(\Fm{}{Z})$ and $\theta_{T'} \in Co_{\bf K}(\Fm{}{Y})$ such that
\begin{enumerate}[(a)]
 \item $T = i(\theta_{T})$, $T' = i(\theta_{T'})$ and $T'' = i(\theta_{T''})$ 
 \item $\theta_{T} \cap \Fm{}{Z} = \theta_{T''} = \theta_{T'} \cap \Fm{}{Z} $
\end{enumerate}
 \end{Df}

The following theorem is the justification for the notion we just introced. 
 
\begin{Teo}\label{ida}\label{TheoremMatrixAmalgamationImpliesCraigInterpolation}
Let $\Sigma$ be a signature, ${\bf K}\subseteq\Sigma\text{-Str}$ a class of algebras closed under subalgebras and $(Co_{{\bf K}},i)$ an equational filter pair having the theory lifting property. If ${\bf K}$ has the amalgamation property, then the logic associated to $(Co_{{\bf K}},i)$ has the Craig entailment interpolation property. 
\end{Teo}

\Dem
By Lemma \ref{2.17} it is enough to prove that $l$ has the flat theory amalgamation property. Let $X,Y$ be non disjoint sets and $T\in Fi_{l_{k}}(\Fm{}{X})$. Denote by $Z=X\cap Y$ and $W=X\cup Y$. Consider $T'=Fi^{Y}_{l_{{\bf K}}}(T\cap \Fm{}{Z})$. So $T'\cap \Fm{}{Z}=T\cap \Fm{}{Z}(=T'')$. Indeed, it is clear that $T\cap \Fm{}{Z}\subseteq T'\cap \Fm{}{Z}$. Suppose $\varphi\not\in T\cap \Fm{}{Z}$, hence $T\cap \Fm{}{Z}\not\vdash_{l_{{\bf K}}}\varphi$. Notice that $Z\subseteq Y$, by \ref{2.19} we have $T\cap \Fm{}{Z}\not\vdash^{Y}\varphi$, then $\varphi\not\in Fi_{l_{{\bf K}}}^{Y}(T\cap \Fm{}{Z})=T'$, proving that if $T'\cap \Fm{}{Z}\subseteq T\cap \Fm{}{Z}$.

From the assumed theory lifting property we obtain $\theta_T \in Co_{\bf K}(\Fm{}{X})$, $\theta_{T''} \in Co_{\bf K}(\Fm{}{Z})$ and $\theta_{T'} \in Co_{\bf K}(\Fm{}{Y})$ such that $T = i(\theta_{T})$, $T' = i(\theta_{T'})$ and $T'' = i(\theta_{T''})$ and $\theta_{T} \cap \Fm{}{Z} = \theta_{T''} = \theta_{T'} \cap \Fm{}{Z} $.

The inclusion $\theta_{T} \cap \Fm{}{Z} \supseteq \theta_{T''}$ tells us that the map $\Fm{}{Z} \hookrightarrow \Fm{}{X}$ induces a well-defined homomorphism $h_1 \colon \Fm{}{Z}/\theta_{T''} \to \Fm{}{X}/\theta_{T}$, and the opposite inclusion $\theta_{T} \cap \Fm{}{Z} \subseteq \theta_{T''}$ tells us that this homomorphism is injective. Likewise, from the equation $\theta_{T''} = \theta_{T'} \cap \Fm{}{Z}$ we obtain a well-defined and injective homomorphism $h_2 \colon \Fm{}{Z}/\theta_{T''} \to \Fm{}{Y}/\theta_{T'}$. 

By hypothesis ${\bf K}$ has the amalgamation property, so we can complete these two maps to a commutative square of injective homomorphisms, for some algebra $\tilde{A} \in {\bf K}$.
\[
\xymatrix{
\Fm{}{Z}/\theta_{T''} \ar@{>->}[d]_{h_1} \ar@{>->}[r]^{h_2} &  \Fm{}{Y}/\theta_{T'} \ar@{>->}[d]^{g_2} \\
 \Fm{}{X}/\theta_{T} \ar@{>->}[r]_{g_1} & \tilde{A}
}
\]  

Since $F \colon Set \to \Sigma\text{-Str}$ is a left adjoint, we have that $\Fm{}{X\cup Y}$ is the pushout of $\Fm{}{X} \hookleftarrow \Fm{}{Z} \hookrightarrow \Fm{}{Y}$. Hence there is a natural homomorphism $\Fm{}{X\cup Y} \to \tilde{A}$. Its image is a subalgebra $A \subseteq \tilde{A}$, hence, by hypothesis, also belongs to ${\bf K}$. Thus $A \cong \Fm{}{X\cup Y}/\theta$ for some congruence $\theta \in Co_{\bf K}(\Fm{}{X\cup Y})$. We claim that the filter on $\Fm{}{X\cup Y}$ required for the flat theory amalgamation property is given by $R:=i_{\Fm{}{X\cup Y}}(\theta)$.

First note that the images of $g_1, g_2$ lie inside $A$, because $\text{Im}\, g_1 = \text{Im}\, (\Fm{}{X} \twoheadrightarrow \Fm{}{X}/\theta_T \to \tilde{A}) = \text{Im}\,(\Fm{}{X} \to \Fm{}{X\cup Y} \to \tilde{A})$ and likewise for $\text{Im}\, g_2$. We can thus in the above square replace $\tilde{A}$ by $A$.

Next, note that for $\varphi \in \Fm{}{X}$ we have $$\varphi \in T\ \ \Leftrightarrow\ \ \delta(\varphi)=\epsilon(\varphi) \in \Fm{}{X}/\theta_T \ \ \Leftrightarrow\ \ \delta(g_1(\varphi))=\epsilon(g_1(\varphi)) \in \Fm{}{X\cup Y}/\theta\ \ \Leftrightarrow\ \  g_1(\varphi) \in R.$$
Here the middle equivalence holds because of $\delta(\varphi)=\epsilon(\varphi)\ \ \Leftrightarrow \ \ g_1(\delta(\varphi))=g_1(\epsilon(\varphi))$ (where $\Leftarrow$ holds because of the injectivity of $g_1$) and $g_1(\delta(\varphi))=\delta(g_1(\varphi))$, $g_1(\epsilon(\varphi))=\epsilon(g_1(\varphi))$ (because $g_1$ is a homomorphism).
Thus $R \cap \Fm{}{X} = T$. With exactly the same reasoning we obtain $R \cap \Fm{}{Y} = T'$. We have proved the flat amalgamation property, and thus the claim.
\qed
 
In view of Theorem \ref{TheoremMatrixAmalgamationImpliesCraigInterpolation} we would like to know when an equational filter pair satisfies the theory lifting property. The following Lemma provides a sufficient criterion for this.

\begin{Lem}\label{LemmaTheoryLiftingHoldsIfXiIsNatural}
Let $(Co_{{\bf K}},i)$ be a congruence filter pair for which the collection of left adjoints $(\Xi_A)_{A \in \Sigma\textrm{-Str}}$ is a natural transformation with respect to variable inclusions, i.e. such that for every inclusion of sets $Z \subseteq X$ the following diagram (in which $j_X \colon \Fm{}{Z} \hookrightarrow \Fm{}{X}$ denotes the induced map of formula algebras) commutes:
\[
\xymatrix{
\Fm{}{Z}\ar[d]^{j_X} & Co_{{\bf K}}(\Fm{}{Z}) && Fi_{l}(\Fm{}{Z}) \ar[ll]^{\Xi_{\Fm{}{Z}}} \\
\Fm{}{X}&Co_{{\bf K}}(\Fm{}{X})\ar[u]^{Co_{\bf K}(j_{X})} && Fi_{l}(\Fm{}{X})\ar[u]_{j^{-1}_{X}} \ar[ll]^{\Xi_{\Fm{}{X}}}
}
\]  
Then $(Co_{{\bf K}},i)$ has the theory lifting property.
\end{Lem}
\Dem 
Suppose we are given sets $X,Y$ such that $Z:=X \cap Y \neq \emptyset$ and a theory $T \subseteq \Fm{}{X}$ of the logic associated to $(Co_{{\bf K}},i)$, define $T'':=T \cap \Fm{}{Z}$ and $T' \subseteq \Fm{}{Y}$ to be the smallest $\Fm{}{Y}$-theory containing $T''$. 
 
First note that $T'\cap \Fm{}{Z}=T\cap \Fm{}{Z}(=T'')$. Indeed, it is clear that $T\cap \Fm{}{Z}\subseteq T'\cap \Fm{}{Z}$. For the other inclusion suppose $\varphi\not\in T\cap \Fm{}{Z}$. Then $T\cap \Fm{}{Z}\not\vdash_{l_{{\bf K}}}\varphi$. As $Z\subseteq Y$, by Prop. \ref{2.19} (which says that the logic based on the set of variables $Y$ is a conservative extension of the logic based on the set of variables $Z$) we have $T\cap \Fm{}{Z}\not\vdash^{Y}\varphi$, hence $\varphi$ is not contained in the closure $T\cap \Fm{}{Z}$, i.e. in $T'$, proving that if $T'\cap \Fm{}{Z}\subseteq T\cap \Fm{}{Z}$.

Denote by $j_X$ (resp. $j_Y$) the map $\Fm{}{Z} \hookrightarrow \Fm{}{X}$ (resp. $\Fm{}{Z} \hookrightarrow \Fm{}{Y}$) induced by the inclusion $Z \subseteq X$ (resp. $Z \subseteq Y$).
Since $i$ is a natural transformation we have the following commutative diagrams.

\[
\xymatrix{
\Fm{}{Z}\ar[d]^{j_X}&Co_{{\bf K}}(\Fm{}{Z})\ar[r]^{i^{{\bf K}}_{Z}}&Fi_{l}(\Fm{}{Z})&&\Fm{}{Z}\ar[d]^{j_Y}&Co_{{\bf K}}(\Fm{}{Z})\ar[r]^{i^{{\bf K}}_{Z}}&Fi_{l}(\Fm{}{Z})\\
\Fm{}{X}&Co_{{\bf K}}(\Fm{}{X})\ar[u]^{Co_{\bf K}(j_{X})}\ar[r]_{i^{{\bf K}}_{X}}&Fi_{l}(\Fm{}{X})\ar[u]_{j^{-1}_{X}}&&\Fm{}{Y}&Co_{{\bf K}}(\Fm{}{Y})\ar[u]^{Co_{\bf K}(j_{Y})}\ar[r]_{i^{{\bf K}}_{Y}}&Fi_{l}(\Fm{}{Y})\ar[u]_{j^{-1}_{Y}}
}
\]  

Note that for any $\theta\in Co_{{\bf K}}(\Fm{}{X})$, we have $Co_{\bf K}(j_{X})(\theta)=\theta\cap \Fm{}{Z}^{2}$ (and the same for $j_{Y}$) and similarly for any $S \subseteq \Fm{}{X}$ we have $j_{X}^{-1}(S)=S\cap \Fm{}{Z}$.

We now construct the congruences of the claim and prove that they satisfy conditions (a) and (b).

By Prop. \ref{FreeAlgebraResult} (which says that for free algebras $i$ is surjective onto filters) there is a congruence $\theta \in Co_{{\bf K}}(\Fm{}{X})$ such that $i_{X}(\theta)=T$. We then know that there is also a minimal such congruence, namely $\theta_{T}:=\Xi_{\Fm{}{X}}(T)$. Define $\theta_{T''}:=\theta_{T}\cap \Fm{}{X}^{2}$.

By the naturality of $i$ we have $i_{\Fm{}{Z}}(\theta_{T''})=i_{\Fm{}{Z}}(\theta_{T}\cap \Fm{}{Z}^{2})=i_{\Fm{}{X}}(\theta_{T})\cap \Fm{}{Z}=T\cap \Fm{}{Z}=T''$. Thus we have constructed $\theta_T \in Co_{\bf K}(\Fm{}{X})$ and $\theta_{T''} \in Co_{\bf K}(\Fm{}{Z})$ and verified that they satisfy conditions (a) and (b) of the claim.

Since $j_{Y}$ is a split monomorphism (because it is induced by an injective map between sets of variables), we have that $Co_{\bf K}(j_{Y})$ is a split epimorphism, hence surjective. Therefore there exists a $\theta'\in Co_{{\bf K}}(\Fm{}{Y})$ such that $\theta'\cap \Fm{}{Z}^{2}=Co_{\bf K}(j_{Y})(\theta')=\theta_{T''}$. Again by naturality of $i$ we have that $T\cap \Fm{}{Z}=T''=i_{Z}(\theta_{T''})=i_{Z}(\theta'\cap \Fm{}{Z}^{2})=i_{Y}(\theta')\cap \Fm{}{Z}$. Thus $T\cap \Fm{}{Z}\subseteq i_{Y}(\theta')$. Since $T'$ is the smallest theory in $\Fm{}{Y}$ containing $T \cap \Fm{}{Z}$, we have $T'\subseteq i_{Y}(\theta')$.

By proposition \ref{FreeAlgebraResult} we have that there exists a $\theta\in Co_{{\bf K}}(\Fm{}{Y})$ such that $i_{Y}(\theta)=T'\subseteq i_{Y}(\theta')$. We now define $\theta_{T'}:=\theta\cap \theta'$. Then $i_{Y}(\theta_{T'})=i_{Y}(\theta)\cap i_{Y}(\theta')=T' \cap i_{Y}(\theta')=T'$, so $\theta_{T'}$ satisfies condition (a).

Next, note that $\theta_{T''}$ is the \emph{minimal} congruence which is mapped to $T''$ by $i_{\Fm{}{Z}}$. Indeed, by assumption $\Xi$ is a natural transformation with respect to homomorphisms induced by variable inclusions. Therefore $\theta_{T''} = \theta_{T}\cap \Fm{}{X}^{2} = Co_{\bf K}(j_{X})(\theta_{T}) = Co_{\bf K}(j_{X})(\Xi(T)) = \Xi(j^{-1}_{X}(T))=\Xi(T'')$. 

Furthermore, we have that $\theta_{T'}\cap \Fm{}{Z}^{2}=(\theta\cap \theta')\cap \Fm{}{Z}^{2}=\theta\cap(\theta'\cap \Fm{}{Z}^{2})=\theta\cap \theta_{T''}$. Thus $\theta_{T'}\cap \Fm{}{Z}^{2}\subseteq \theta_{T''}$. 
We also have $i_Y(\theta_{T'}\cap \Fm{}{Z}^{2}) = i_Y(Co_{\bf K}(j_Y)(\theta_{T'})) = j_Y^{-1}(i_Y(\theta_{T'})) = j_Y^{-1}(T') = T' \cap \Fm{}{Z} = T''$, and since $\theta_{T''}$ was the \emph{minimal} congruence which is mapped to $T''$ by $i_Z$ this implies $\theta_{T''} \subseteq \theta_{T'}\cap \Fm{}{Z}^{2}$. Altogether we have $\theta_{T'}\cap \Fm{}{Z}^{2} = \theta_{T''}$ and thus $\theta_{T'}$ satisfies condition (b).
\qed

\begin{Obs}
 The naturality condition of Lemma \ref{LemmaTheoryLiftingHoldsIfXiIsNatural} can be spelled out as the requirement that \linebreak $\Xi_{\Fm{}{X}}(j_X^{-1}(T)) \subseteq (j_X \times j_X)^{-1}(\Xi_{\Fm{}{Y}}(T))$.
 By Lemma \ref{LemmaSubstitutionInvarianceForTheEquationalOperators}(\ref{substitutionInvarianceInverseImageXi}), the inclusion $\Xi_{\Fm{}{X}}(\sigma^{-1}(T)) \subseteq (\sigma\times \sigma)^{-1}(\Xi_{\Fm{}{Y}}(T))$ always holds and one only needs to check the other inclusion.
\end{Obs}

\begin{Cor}\label{CorollaryiInjectiveImpliesTheoryLifting}
 Let $(Co_{{\bf K}},i)$ be a congruence filter pair with $i$ injective. Then $(Co_{{\bf K}},i)$ has the theory lifting property.
\end{Cor}
\Dem
  By Prop. \ref{CriterionAlgebraizable} the logic presented by $(Co_{{\bf K}},i)$ is algebraizable with $i$ being a substitution preserving lattice isomorphism. By Theorem \ref{TheoremSecondIsoThmInclUniquenessOfTheQuasivariety} the inverse is given by $\Omega$, which hence is a left adjoint. By uniqueness of adjoints we have $\Xi_A = \Omega_A$ for every $\Sigma$-structure $A$ (see also Theorem \ref{TheoremInclusionsXiIdLeibniz}). Since the logic associated to $(Co_{{\bf K}},i)$ is algebraizable, it is in particular equivalential and hence by \cite[Thm 6.68]{Fon} $\Omega$ is a natural transformation. Hence $\Xi = \Omega$ satisfies the condition of Lemma \ref{LemmaTheoryLiftingHoldsIfXiIsNatural}.
\qed 

We arrive at the following well-known fact about algebraizable logics.

\begin{Cor}\label{CorollaryForAlgebraizableLogicsAmalgamationImpliesCraigInterpolation}
 Let $l$ be an algebraizable logic with equivalent semantics in a quasivariety ${\bf K}$. If ${\bf K}$ has the amalgamation property, then the $l$ has the Craig entailment property. 
\end{Cor}
\Dem
Since the logic $l$ is algebraizable, it can be presented by an equational filter pair with $i$ injective, just consider the pair $(Co_{K}, i)$ such that K is the equivalent algebraic semantic for $l$ and $i_{A}$ is the inverse of $\Omega^{A}$ for any algebra $A$. Now one can apply Corollary \ref{CorollaryiInjectiveImpliesTheoryLifting} and Theorem \ref{TheoremMatrixAmalgamationImpliesCraigInterpolation}.
\qed 

Of course stronger statements than the preceding one are known, e.g. that amalgamation is equivalent to Craig interpolation for algebraizable logics, see \cite[Cor. 5.27]{CzP}. 
While we could obtain Corollary \ref{CorollaryForAlgebraizableLogicsAmalgamationImpliesCraigInterpolation} as a quick byproduct of Lemma \ref{LemmaTheoryLiftingHoldsIfXiIsNatural}, the aim of the latter is to pave the way for new cases.

We now give another criterion for the naturality condition of Lemma \ref{LemmaTheoryLiftingHoldsIfXiIsNatural} to be satisfied. For this recall from Section \ref{SubsectionCongruencesAndEquationalConsequence} the notation $Cn_{\bf K}^X (\Gamma)$ for the ${\bf K}$-congruence relation generated by a set of equations $\Gamma$ on the free algebra $\Fm{}{X}$ or, equivalently, for the closure operator for the equational consequence relation. 
 
\begin{Lem}\label{LemmaConditionForNaturalityOfXi}
 Let ${\bf K}$ be a quasivariety and $\langle \delta(\varphi), \epsilon(\varphi) \rangle$ an equation over its signature. Suppose that every inclusion $\sigma \colon \Fm{}{Z} \hookrightarrow \Fm{}{X}$ of free algebras coming from an inclusion of generators $Z \subseteq X$ induces a conservative extension of equational theories, i.e. that
   for every congruence $\theta \in Co_{\bf K}(\Fm{}{X})$ we have  
 $$(Cn_{\bf K}^X\{ \langle \delta(\varphi), \epsilon(\varphi) \rangle \in \theta \mid \varphi \in \Fm{}{X} \})  \cap \Fm{}{Z}^2 
  = Cn_{\bf K}^Z(\{ \langle \delta(\varphi), \epsilon(\varphi) \rangle \in \theta \mid \varphi \in \Fm{}{Z} \}$$
   
 Then the equational filter pair associated to the equation satisfies the naturality condition of Lemma \ref{LemmaTheoryLiftingHoldsIfXiIsNatural}, i.e. $(\sigma \times \sigma)^{-1}(\Xi_{\Fm{}{X}}(T)) = \Xi_{\Fm{}{Z}}(\sigma^{-1}T)$
\end{Lem}
\Dem
First note that $\sigma^{-1}$ (resp. $(\sigma \times \sigma)^{-1}$) is just restriction to formulas (resp. pairs of formulas) in $\Fm{}{Z}$, i.e. $\sigma^{-1}(T) = T \cap \Fm{}{Z}$ and $(\sigma \times \sigma)^{-1}(\theta) = \theta \cap \Fm{}{Z}^2$. Next note that $\{ \langle \delta(\varphi), \epsilon(\varphi) \rangle \in \theta \mid \varphi \in \Fm{}{X} \} = \{ \langle \delta(\varphi), \epsilon(\varphi) \rangle \mid \varphi \in i(\theta) \}$

As by Prop. \ref{FreeAlgebraResult} for free algebras filters are exactly the sets of formulas of the form $T=i(\theta)$ for some $\theta \in Co_{\bf K}(\Fm{}{X})$, we can rewrite the assumption of conservative extension of equational theories by saying that for every filter $T$ we have
$$(Cn_{\bf K}^X\{ \langle \delta(\varphi), \epsilon(\varphi) \rangle \mid \varphi \in T \})  \cap \Fm{}{Z}^2 = Cn_{\bf K}^Z(\{ \langle \delta(\varphi), \epsilon(\varphi) \rangle \mid \varphi \in T \}  \cap \Fm{}{Z}^2) \ \ \ \ \ (\ast \ast)$$ 

Now the claim follows from the series of equations
$$\begin{array}{rclr}
(\sigma \times \sigma)^{-1}(\Xi_{\Fm{}{X}}(T)) &=& \Xi_{\Fm{}{X}}(T) \cap \Fm{}{Z}^2 & \\
  &=& (Cn_{\bf K}^X\{ \langle \delta(\varphi), \epsilon(\varphi) \rangle \mid \varphi \in T \})  \cap \Fm{}{Z}^2 & \text{(Cor. \ref{PropOurLeftAdjointIsBlokPigozzisOmegaK})}\\
  &=& Cn_{\bf K}^Z(\{ \langle \delta(\varphi), \epsilon(\varphi) \rangle \mid \varphi \in T \}  \cap \Fm{}{Z}^2) & (\ast \ast) \\
  &=& Cn_{\bf K}^Z(\{ \langle \delta(\varphi), \epsilon(\varphi) \rangle \mid \varphi \in T \cap \Fm{}{Z}\}  ) &  \\
  &=& \Xi_{\Fm{}{Z}}(T \cap \Fm{}{Z}) & \text{(Cor. \ref{PropOurLeftAdjointIsBlokPigozzisOmegaK})}\\
  &=& \Xi_{\Fm{}{Z}}(\sigma^{-1}T)
\end{array}$$
\qed

Note that the condition of Lemma \ref{LemmaConditionForNaturalityOfXi} is purely a condition on the quasivariety ${\bf K}$ and the equations, involving no logic or particularities of filter pairs.

%
%

\begin{Ex}\label{ExampleSuccessorLogicHasTheoryLiftingProperty}
Consider again the logic $l_s$ of Example \ref{ExampleSuccessoralgebra} associated to the equational filter pair $(Co_{\bf K}, i^\tau)$ where $\Sigma$ is the signature with just one unary function symbol $s$, ${\bf K}$ is the class of all $\Sigma$-structures and $\tau$ is the single equation $\langle x, s(x) \rangle$. We will show that the conditions of Lemma \ref{LemmaConditionForNaturalityOfXi} are satisfied.

To this end note that the set $\{ \langle \varphi, s(\varphi) \rangle \in \theta \mid \varphi \in \Fm{}{X} \}$ occurring in this condition is the set of $\varphi \in \Fm{}{X}$ which are mapped to $s$-fixed points in $\Fm{}{X} / \theta$. 
The congruence relation generated by such a set is simply the congruence which identifies all sucessors of such a fixed point. Explicitly $\varphi$ and $\psi$ get identified if and only if they satisfy what one might call the successor-fixed point condition, namely that either $\varphi\textrm{ is an \emph{s}-fixed point in }\Fm{}{X}/\theta\textrm{ and  }\psi=s^{\circ n}(\varphi)$ or $\psi\textrm{ is an \emph{s}-fixed point in }\Fm{}{X}/\theta\textrm{ and  }\varphi=s^{\circ n}(\psi)$ for some $n \in \mathbb{N}$ (here $s^{\circ n}$ denotes the $n$-fold iteration of $s$). The left one of the two congruences occurring in the condition of Lemma \ref{LemmaConditionForNaturalityOfXi} then consists of those pairs $\langle \varphi, \psi \rangle$ which satisfy the said condition and whose variables lie in $Z$.

The corresponding description holds of course also for the right hand side; it is the set of pairs $\langle \varphi, \psi \rangle$ whose variables lie in $Z$ and which satisfy the successor-fixed point condition. Thus the two sets are equal and the hypotheses of Lemma \ref{LemmaConditionForNaturalityOfXi} are satisfied.
\end{Ex}

We thus have a first example of a non-protoalgebraic logic for which we can infer the Craig entailment property from the amalgamation property:

\begin{Ex}\label{ExampleSucessorLogicHasMatrixAmalgamationProperty}
Consider the logic $l_s$ of example \ref{ExampleSuccessoralgebra}. We continue using the notation of loc. cit., i.e. the signature $\Sigma$ is the one with exactly one unary operation and $i$ is the natural transformation from the filter pair defining that logic. The category of $\Sigma$-algebras (i.e. of sets with an endomorphism) is equivalent to the functor category $Set^\mathbb{N}$, where one sees the monoid $(\mathbb{N},+)$ as a category with one object. This is a functor category, thus pushouts are formed on the underlying sets. In particular, pushouts of monomorphisms are monomorphisms, because this is the case in the category of sets. This shows that $\Sigma$-Alg has the amalgamation property.

Example \ref{ExampleSuccessorLogicHasTheoryLiftingProperty} together with Lemma \ref{LemmaConditionForNaturalityOfXi} and Lemma \ref{LemmaTheoryLiftingHoldsIfXiIsNatural} implies that the filter pair we used to define $l_s$ has the theory lifting property.

Thus Theorem \ref{TheoremMatrixAmalgamationImpliesCraigInterpolation} applies and we can deduce that $l_s$ has the Craig entailment property. We have successfully applied algebraic reasoning to an (admittedly artificial) logic that is neither protoalgebraic nor truth-equational, see Examples \ref{ExampleSuccessoralgebra} and \ref{ExampleSuccessorLogicIsNotTruthEquational}.
\end{Ex}

We conclude this section with determining a class of filter pairs for which the theory lifting property holds.

%

\begin{Df}\label{DefregularEquationsAndVarieties}
 An equation of formulas $\langle \varphi, \psi \rangle \in \Fm{}{X}$ is called \emph{regular} if  $\varphi, \psi $ contain exactly the same variables. A variety is called regular if it can be defined by regular equations.
\end{Df}

\begin{Ex}
 The variety of bounded semilattices is regular: it can be axiomatized over the signature $\Sigma=\{ \wedge, \top \}$ by the equations $x \wedge x = x$, $x \wedge y = y \wedge x$, $x \wedge (y \wedge z) = (x \wedge y) \wedge z$ and $x \wedge \top = x$.

 \medskip
 
Some varieties of bounded semilattices with a rudimentary implication operation are regular, e.g the one obtained by enhancing the above signature to $\Sigma'=\{ \to, \wedge, \top \}$ and adding the axioms $x \wedge (x \to y) = x \wedge y$ and $x \to (x \wedge y) = (x \to x) \wedge (x \to y)$. Further typical axioms like $x \to x = \top$ and $y \wedge (x \to y)$ are, however, not regular.
  
 \medskip
 
 The variety of lattices is not regular. Its standard axiomatization contains the absorption law $x \wedge (y \vee x) = x$ which is not regular, in particular the free lattice over some set of generators satisfies that law. If there was an axiomatization by regular quasiequations then, by arguments similar to the ones in the proof of Prop. \ref{PropregularVarietiesHaveTheoryLiftingProperty} below, one could show that the free lattice does not satisfy this equation.
\end{Ex}

\begin{Prop}\label{PropregularVarietiesHaveTheoryLiftingProperty} Let ${\bf K}$ be a regular variety. Then any equational filter pair $(Co_{\bf K}, i^\tau)$, with $\tau = \langle \delta(x), \epsilon(x) \rangle$ and $\delta(x)$ and $\epsilon(x)$ each having exactly the free variable $x$, has the theory lifting property.
\end{Prop}
\Dem
 We verify the hypothesis of Lemma \ref{LemmaConditionForNaturalityOfXi}, i.e.  $$(Cn_{\bf K}^X\{ \langle \delta(\varphi), \epsilon(\varphi) \rangle \in \theta \mid \varphi \in \Fm{}{X} \})  \cap \Fm{}{Z}^2 
  = Cn_{\bf K}^Z(\{ \langle \delta(\varphi), \epsilon(\varphi) \rangle \in \theta \mid \varphi \in \Fm{}{Z} \}.$$
  
  \medskip
  
\noindent  $\supseteq$: We have 
    $$\begin{array}{rclr}
(Cn_{\bf K}^X\{ \langle \delta(\varphi), \epsilon(\varphi) \rangle \in \theta \mid \varphi \in \Fm{}{X} \})  \cap \Fm{}{Z}^2 &\supseteq & (Cn_{\bf K}^X\{ \langle \delta(\varphi), \epsilon(\varphi) \rangle \in \theta \mid \varphi \in \Fm{}{Z} \})  \cap \Fm{}{Z}^2 \\
  &=& Cn_{\bf K}^Z(\{ \langle \delta(\varphi), \epsilon(\varphi) \rangle \in \theta \mid \varphi \in \Fm{}{Z} \} 
  \end{array}$$
where the first inclusion holds because of $\Fm{}{X} \supseteq \Fm{}{Z}$ and the second equality follows from Lemma \ref{LemmaInclusionOfVariablesGivesConservativeTranslationOfEquationalConsequences}.
  
  \medskip
  
\noindent  $\subseteq$: Let $\langle \rho, \psi \rangle \in \Fm{}{Z}^2$ with $\{ \langle \delta(\varphi), \epsilon(\varphi) \rangle \in \theta \} \vDash^X_{\bf K}  \langle \rho, \psi \rangle$. Then by \cite[Thm 2.1]{Cze3} there exists a finite sequence $ \langle \alpha_1, \beta_1 \rangle, \ldots, \langle \alpha_n, \beta_n \rangle$  with $ \langle \alpha_n, \beta_n \rangle = \langle \rho, \psi \rangle$ and such that for every $i = 1,  \ldots , n$ the equation $\langle \alpha_i, \beta_i \rangle$ satisfies one of the following:
\begin{enumerate}
\item one has $\alpha_i = \beta_i$.

\item one has $\langle \alpha_i, \beta_i \rangle = \langle \delta(\varphi), \epsilon(\varphi) \rangle \in \theta$ for some $\varphi \in \Fm{}{X}$

\item there is a substitution mapping one of the defining equations for ${\bf K}$ to $\langle \alpha_i, \beta_i \rangle$.


\item the equation $\langle \alpha_i, \beta_i \rangle$ follows from $\{\langle \alpha_1, \beta_1 \rangle, \ldots, \langle \alpha_{i-1}, \beta_{i-1} \rangle\}$ by one of Birkhoff's rules $\langle \gamma_1, \gamma_2 \rangle \vDash \langle \gamma_2, \gamma_1 \rangle$ (symmetry), $\{ \langle \gamma_1, \gamma_2 \rangle, \langle \gamma_2, \gamma_3 \rangle \} \vDash \langle \gamma_1, \gamma_3 \rangle$ (transitivity) or \linebreak $\{ \langle \gamma_1, \gamma_1' \rangle, \ldots, \langle \gamma_k, \gamma_k' \rangle \} \vDash \langle f(\gamma_1, \ldots, \gamma_k), f(\gamma_1', \ldots, \gamma_k') \rangle$ for a $k$-ary function symbol $f$ (congruence).
\end{enumerate}
First, observe that by induction one can see that every equation $\langle \alpha_i, \beta_i \rangle$ is regular. Indeed, if $\langle \alpha_i, \beta_i \rangle$ is introduced by one of the first three rules, then it is regular because the image of a regular equation under a substitution is regular, and the equations $\langle x, x \rangle$, $\langle \delta(x), \epsilon(x) \rangle$ and the defining equations of ${\bf K}$ are regular, in the first case obviously so, in the second and third case by assumption on ${\bf K}$.

If $\langle \alpha_i, \beta_i \rangle$ is introduced by the rules of the last item, then observe that these transform regular equations, again into regular equations. Since by induction hypothesis $\{\langle \alpha_1, \beta_1 \rangle, \ldots, \langle \alpha_{i-1}, \beta_{i-1} \rangle\}$ are regular, so is $\langle \alpha_i, \beta_i \rangle$.

Now, again by induction on $n$, we show that one can derive $\langle \rho, \psi \rangle$ already from the hypotheses in $\{ \langle \delta(\varphi), \epsilon(\varphi) \rangle \in \theta \mid \varphi \in \Fm{}{Z} \}$, i.e. without using formulas in $\Fm{}{X} \setminus \Fm{}{Z}$. Indeed, if one arrives at the last step $ \langle \alpha_n, \beta_n \rangle = \langle \rho, \psi \rangle$ by application of one of the first three rules, this is clear. 

If one arrives at the last equation by application of the symmetry rule $\langle \psi , \rho \rangle \vDash \langle \rho, \psi \rangle$, then the equation $\langle \psi , \rho \rangle$ also has variables in $Z$ and has a shorter proof, hence by induction hypothesis follows from $\{ \langle \delta(\varphi), \epsilon(\varphi) \rangle \in \theta \mid \varphi \in \Fm{}{Z} \}$.

If one arrives at the last equation by application of the congruence rule \linebreak $\{ \langle \gamma_1, \gamma_1' \rangle, \ldots, \langle \gamma_k, \gamma_k' \rangle \} \vDash \langle f(\gamma_1, \ldots, \gamma_k), f(\gamma_1', \ldots, \gamma_k') \rangle$ for a $k$-ary function symbol $f$, then, since the outcome only contains variables in $Z$, the equations $\langle \gamma_i, \gamma_i' \rangle$ also contain only variables from $Z$. By induction hypothesis all of these equations can be derived from $\{ \langle \delta(\varphi), \epsilon(\varphi) \rangle \in \theta \mid \varphi \in \Fm{}{Z} \}$.

Finally, if one arrives at the last equation by application of the transitivity rule $\{ \langle \rho, \gamma \rangle, \langle \gamma, \psi \rangle \} \vDash \langle \rho, \psi \rangle$, then, since $\rho$ and $\psi$ have variables in $Z$ and $\langle \rho, \gamma \rangle, \langle \gamma, \psi \rangle$ are regular, they also have variables in $Z$, hence by induction hypothesis these equations can be derived  from $\{ \langle \delta(\varphi), \epsilon(\varphi) \rangle \in \theta \mid \varphi \in \Fm{}{Z} \}$.
\qed

\begin{Cor}\label{CorollaryAmalgamationInregularVarietiesImpliesInterpolation}
 Let ${\bf K}$ be a regular variety satisfying the amalgamation property.
 Then any logic $l$ with an algebraic semantics in ${\bf K}$ given by a regular equation satisfies the Craig interpolation property.
\end{Cor}
\Dem
By Theorem \ref{TeoEveryLogicWithAlgebraicSemanticsComesFromEqunlFilterPair} the logic $l$ can be presented by an equational filter pair. By the construction of that filter pair, it will satisfy the hypotheses of Prop. \ref{PropregularVarietiesHaveTheoryLiftingProperty}, so that we have the theory lifting property. By Theorem \ref{TheoremMatrixAmalgamationImpliesCraigInterpolation} it follows that amalgamation implies Craig interpolation.
\qed

We point out that Examples \ref{ExampleSuccessorLogicHasTheoryLiftingProperty}, \ref{ExampleSucessorLogicHasMatrixAmalgamationProperty} are an instance of Corollary \ref{CorollaryAmalgamationInregularVarietiesImpliesInterpolation}.

While the class of logics covered by Corollary \ref{CorollaryAmalgamationInregularVarietiesImpliesInterpolation} does not commonly show up in applications, the result nevertheless shows that systematic connections between the amalgamation property and the Craig interpolation property exist even beyond the classes of protoalgebraic or truth-equational logics. It seems worthwhile to look for further practical criteria for a filter pair to satisfy the theory lifting property, or to satisfy the condition on congruence generation appearing in Lemma \ref{LemmaConditionForNaturalityOfXi}.

\section{Final remarks and future works}\label{SectionVista}

We sketch a few directions into which the work on filter pairs may be taken. We begin with three sets of questions upon which we have touched in this work and which may be developed further.


\begin{Obs}\label{Tomasso Jansana } {\bf Adjoint operators and the Leibniz hierarchy}
 With the sample of results of section \ref{filteradjunctionequationalfilterpairs}, the topic of adjunctions and filter pairs is far from exhausted. While we have considered the inclusions $\theta \subseteq \Omega(i(\theta))$, $i(\theta) \subseteq i(\Omega(i(\theta)))$ etc. and the meaning of equality holding in such an inclusion, one can, for example, also consider when equality holds for the opposite compositions $i(\Omega(F))=F$, $\Omega (i(\Omega(F)))=\Omega (F)$ and so on. Concrete descriptions of the operator $\Omega \circ i$ and its companions can be given for equational filter pairs, with assertional logics being the easiest case. Each of these conditions will define a class of logics, possibly extending the Leibniz hierarchy.
 
 A natural further step will be to give corresponding treatments of the Suszko operator and work out the relation to the strong version of a logic.

 Furthermore, in \cite{JM1}, \cite{JM2} there is an interesting approach to the Leibniz hierarchy: Roughly, Leibniz classes are taken to be classes defined by the condition of receiving a translation from a fixed class of logics. It is an interesting question whether the classes of logics defined via properties of the operators $i, \Omega, \Xi$ and their adjoints, are Leibniz classes in the sense of \cite{JM1},  \cite{JM2}.
\end{Obs}

\begin{Obs}\label{Rem bridge-theorems} {\bf Bridge theorems through equational filter-pairs.}
In Section \ref{SectionCraigInterpolation} we have started to look into results connecting amalgamation and Craig interpolation. One may try to refine these results, for example find new criteria for the theory lifting property to hold. The theory lifting does not need to occur through the operators $\Xi$ or $\Omega$; these two are just the extreme ends between which further operators might be found with the necessary naturality condition. Even concentrating on the operator $\Xi$ there is much room for further investigations. The question when the condition of Lemma \ref{LemmaConditionForNaturalityOfXi} holds is purely a matter of universal algebra, that is worth investigating. For particular quasivarieties and equations an analysis of the proofs in the associated equational logic, along the lines of the one carried out in the proof of Proposition \ref{PropregularVarietiesHaveTheoryLiftingProperty}, may lead to amalgamation results.

It will also be interesting to look for further results connecting algebraic properties with logical properties for congruence filter pairs, e.g. the Beth property, or being admissibly closed, which for an algebraizable logic can be expressed in terms of the category of matrices or the category found in \cite{MaPi2}.
\end{Obs}

\begin{Obs}\label{Rem Logic by equations}{\bf Studying logics by variation of their defining equation.}
A common instance of the situation of Theorem \ref{TheoremLogicsFromEquations} is when the signature $\Sigma$ contains connectives $\wedge, \vee, \top, \bot, \neg$ for conjunction, disjunction, truth, falsity and negation (and possibly others), ${\bf K}$ is a quasivariety of lattices with extra operations and the equation of the filter pair is simply given by $\langle x, \top \rangle$. This means that filters are sets of elements which are the equivalence class of the element $\top$ in some quotient.

Even if one is only interested in this logic, one may get insight by varying the equation. For example, maintaining the same variety ${\bf K}$, one has also the logics associated to the equations $\langle x, \neg x \rangle$, $\langle x \wedge \neg x, \bot \rangle$, $\langle x \vee \neg x, \top \rangle$ or $\langle x , \neg\neg x \rangle$. The theorems of these logics can be seen to be the ``contingent formulas'' (which are ``as true as they are false'', which e.g. makes sense in many-valued logic), the consistent formulas (which is a non-trivial property for paraconsistent logics), the exhaustive formulas (which, together with their negation exhaust all possible cases) or the classical formulas (e.g. considered in intuitionistic logic as the image of the double negation translation).

The consequence relation of these new logics can, roughly, be understood as $\Gamma \vdash_{cons} \varphi$ saying that if all $\gamma \in \Gamma$ are, say, consistent, then so is $\varphi$. The relationship to the initially given logic is ensured by taking the same quasivariety ${\bf K}$.

One can start exploring the new associated logics with the help of the congruence filter pair presentations. For example, suppose that one knows that the initial logic is algebraizable and enjoys the Craig interpolation property. Then one knows that ${\bf K}$ has the amalgamation property and can try to infer Craig interpolation for the new logics by the techniques of Section \ref{SectionCraigInterpolation}.
\end{Obs}

While the previous directions of further research are close to the contents of this article, there are also quite different directions into which the idea of filter pairs may be usefully developed.

\begin{Obs}\label{Rem new hierarchy for filter pairs} {\bf Further examples of filter pairs.}
In this work we have emphasized the class of congruence filter pairs. But there are other interesting classes. For example one can consider, instead of the lattice of (relative) congruences on algebras, the lattice of (relative) congruences for some subsignature (or more generally a morphism of signatures with target the signature of the given logic).

\medskip

Let us consider the extreme example of the empty subsignature: Then the functor part of the filter pair will associate to an algebra the lattice of equivalence relations on its underlying set. In general a quotient by such an equivalence relation does not inherit an algebra structure. It \emph{does}, however, inherit a structure of \emph{multialgebra}, i.e. a gadget with many-valued operations corresponding to the connectives of the initially given signature. This could connect to Avron's non-deterministic semantics \cite{Avron} \cite{AvronZamanskyHandbook}, and maybe make it possible to relate metalogical properties to non-deterministic matrices in the style of algebraic logic.

\medskip

Going further, one may take a subsignature and consider consequence relations only with respect to that subsignature. Repeating the considerations of the previous section, one could hope to obtain a partial Craig interpolation theorem, talking only about the fragment given by the subsignature.

\medskip

Finally, allowing more general signature morphisms, and considering congruences for the pullbacks of algebras along these, one may obtain a framework for remote algebraization from the left -- the dual to the possible translation semantics of \cite{BCC1} \cite{BCC2}.

\medskip

Also any functorial reflective sublattice of the lattice of all congruences can be considered, as long as the inclusion preserves directed suprema. The passage to a reflective sublattice will result in a stronger logic (to see this just consider the composition of the two adjunctions in question).


\medskip

Of course, lattices do not need to consist of congruences at all. One possible direction, still staying on the side of algebra, is to let ${\bf K}$ be any reflexive subcategory of $\Sigma$-Str, not necessarily a quasivariety, and for a $\Sigma$-algebra $A$ consider a lattice of maps into certain members of ${\bf K}$ ordered by factorization. In a different direction, both the quotient and the subobject lattices in a locally finitely presentable category are algebraic lattices \cite[Thm. 5, Thm. 12]{PorstAlgLattices}. One might try to take the subobject lattices in a category of coalgebras and obtain a new point of view on coalgebraic semantics. 
\end{Obs}

\begin{Obs}
{\bf Algebraic semantics in other categories than Sets.}
One may adjust the notion of filter pair by changing the domain category from $\Sigma$-structures in the category of sets to $\Sigma$-structures in different categories. Not only will this extend the range of semantics that can be given to propositional logics. It should also render the statement of Prop. \ref{CorollaryAmalgamationInregularVarietiesImpliesInterpolation} far more useful: While regular varieties defined over the category of sets are scarce, it is precisely regularity that allows to interpret such varieties in arbitrary symmetric monoidal categories. For example, the variety of rings is not regular, but the variety of monoids is. Hence one can say what a monoid in a symmetric monoidal category is. Taking the category of abelian groups with tensor product, one recovers the category of rings.
\end{Obs}






\begin{Obs}\label{Rem InfFilterPair+Category} {\bf On categories of infinitary logics and infinitary filter pairs.}
Motivated by a question on the existence of natural extensions of a logic, posed by Cintula and Noguera in \cite{CN} and in the meantime solved by P{\v r}enosil \cite{Prenosil}, we have extended in \cite{AMP2} the notion of finitary filter pair to an infinitary version called $\kappa$-filter pair, where $\kappa$ is any regular cardinal. We show that a $\kappa$-filter pair gives rise to a logic of cardinality $\kappa$ and every logic of cardinality $\kappa$ comes from a $\kappa$-filter pair, extending the result for $\kappa = \omega$ contained in this work. Together with an analogue of Prop. \ref{PropHomomorphismsInduceTranslationsAndVariableInclusionsInduceConservativeTranslations} this yields a solution to the construction of natural extensions, which turns out to be identical to one of P{\v r}enosil's solutions, but obtained on a different route. 

Due to the possible non-uniqueness of natural extensions, the ambiguity of different presentations of a single logic by different filter pairs is bigger than in the finitary case. This suggests a point of view of a filter pair as a logic together with a family of natural extensions to sets of variables of all arities.

The natural task arises to extend the results of this work to the infinitary case. This is carried out in \cite{AMP2} for the results of Section \ref{SectionFilterFunctors}.
\end{Obs}

\begin{Obs}\label{Rem MetaLogic} {\bf On a representation theory of logics.}
Taking adequate notions of morphisms, we show that the category of logics of cardinality $\kappa$ and translation morphisms is (isomorphic to) a full reflective subcategory of the category of $\kappa$-filter pairs: Thus this is can be see as a natural sequel of the work initiated in \cite{MaPi2}, where there is established a full correspondence of certain functors between categories of $\Sigma$-structures and translations between algebraizable logics. 

The presentation of functorial encodings of morphisms of logics,  will play a role in the long term project of on establishing a representation theory of logics, i.e. studying arbitrary logics through their translations into a class of "well-behaved logical objects", such as the class of algebraizable logics (\cite{MaPi1}, \cite{MaPi2}) or  the class of  filter pair matrices (\cite{AMP1}).

A possibly interesting use of these machinery is to establish local-global principles in the realm of Propositional Logic, in particular in the study of meta-logical properties of logical systems. A concrete goal is to establish sufficient conditions for a meta-logical property to be preserved under constructions such as products, filtered colimits, among others. 
\end{Obs}

\end{document}